\documentclass[a4paper,10pt,draft]{article}
\usepackage{amssymb, amsmath, amsthm, latexsym, verbatim}

\newcommand \n{\nabla}
\newcommand \la{\lambda}

\newcommand \Rn{\mathbb R^n}
\newcommand {\Ro}{\mathbb R^8}

\newcommand \rk{\mathrm{rk}}
\newcommand \Ker{\mathrm{Ker}}

\newcommand \Cliff{\mathrm{Cliff}}
\newcommand \Cl{\mathrm{Cl}}
\newcommand \Span{\mathrm{Span}}
\newcommand \Tr{\mathrm{Tr}}

\newcommand \diag{\mathrm{diag}}

\newcommand \T{\theta}
\renewcommand \a{\alpha}

\newcommand \<{\langle}
\renewcommand \>{\rangle}

\theoremstyle{plain}

\newtheorem*{theorem*}{Theorem}
\newtheorem*{corollary*}{Corollary}
\newtheorem{lemma}{Lemma}
\newtheorem{proposition}{Proposition}
\newtheorem{fact}{Fact}
\newtheorem*{namedtheorem}{\theoremname}
\newcommand{\theoremname}{te}
\newenvironment{named}[1]{\renewcommand{\theoremname}{#1}
\begin{namedtheorem}}{\end{namedtheorem}}

\theoremstyle{definition}
\newtheorem{definition}{Definition}

\theoremstyle{remark}
\newtheorem*{remark}{Remark}

\begin{document}

\title{Osserman manifolds of dimension $8$}

\author{Y.Nikolayevsky\thanks{The work is supported by ARC Discovery grant DP0342758}}

\date{}


\maketitle

\begin{abstract}
For a Riemannian manifold $M^n$ with the curvature tensor $R$, the
Jacobi operator $R_X$ is defined by $R_XY = R(X,Y)X$.
The manifold $M^n$ is called {\it pointwise Osserman} if, for every
$p \in M^n$, the eigenvalues of the Jacobi operator $R_X$  do not depend of a
unit vector $X \in T_pM^n$, and is called {\it globally Osserman} if they
do not depend of the point $p$ either. R.~Osserman conjectured that globally
Osserman manifolds are flat or rank-one symmetric. This Conjecture is true for
manifolds of dimension $n \ne 8, 16$ \cite{NarXiv}.
Here we prove the Osserman Conjecture and its
pointwise version for $8$-dimensional manifolds.
\end{abstract}

\section{Introduction}
\label{intro}

An {\it algebraic curvature tensor} $R$ in a Euclidean space $\Rn$ is a $(3, 1)$
tensor having the same symmetries as the curvature tensor of a Riemannian manifold.
For $X \in \Rn$, the {\it Jacobi operator} $R_X : \Rn \to \Rn$ is defined by
$R_XY = R(X, Y)X$ . The Jacobi operator is symmetric and $R_XX = 0$ for all
$X \in \Rn$. Throughout the paper, ``eigenvalues of the Jacobi operator"
refers to eigenvalues of the restriction of $R_X$, with $X$ a
unit vector, to the subspace $X^\perp$.

\begin{definition} An algebraic curvature tensor $R$ is called
{\it Osserman} if the eigenvalues of the Jacobi operator $R_X$ do not
depend of the choice of a unit vector $X \in \Rn$.
\end{definition}

\begin{definition} A Riemannian manifold $M^n$ is called
{\it pointwise Osserman} if its curvature tensor is Osserman. If, in
addition, the eigenvalues of the Jacobi operator are constant on $M^n$,
the manifold $M^n$ is called {\it globally Osserman}.
\end{definition}

Flat and rank-one symmetric spaces
are globally Osserman, since the isometry group of each of them
acts transitively on its unit tangent bundle. Osserman \cite{O}
conjectured that the converse is also true:

\begin{named}{Osserman Conjecture}
A globally Osserman manifold is flat or rank-one symmetric.
\end{named}

In the most cases, the answer to the Osserman Conjecture is affirmative, as well as
to its ``pointwise" version (see the Corollary below).

In this paper, we prove the Osserman Conjecture for $8$-dimensional manifolds:

\begin{theorem*}
A pointwise Osserman manifold of dimension eight is flat or rank-one symmetric.
\end{theorem*}

Combining this with Theorems 1 and 2 of \cite{NarXiv}, we get:

\begin{corollary*} In each of the following cases a Riemannian manifold $M^n$ is flat or
rank-one symmetric:

1) $M^n$ is globally Osserman and $n \ne 16$.

2) $M^n$ is pointwise Osserman and $n \ne 2, 4, 16$.

3) $n=16$, the manifold $M^{16}$ is (pointwise or globally) Osserman, and
its Jacobi operator has no eigenvalues of multiplicity $m \in \{7, 8, 9\}$.
\end{corollary*}

In the cases covered by the Corollary, there is not much difference between globally and pointwise
Osserman conditions, except in dimension $2$, where any Riemannian manifold is pointwise Osserman,
and in dimension $4$, where any globally Osserman manifold is flat or rank-one symmetric \cite{Chi},
but there exist pointwise Osserman manifolds that are not symmetric
(``generalized complex space forms", see \cite[Corollary 2.7]{GSV}, \cite{Ols}).

We refer to \cite{GKV} for results on the Osserman Conjecture in Riemannian and
semi-Riemannian geometry.

The paper is organized as follows. In Section~\ref{sec:2}, we give the proof
of the Theorem, assuming two Propositions: Proposition~\ref{p1} in which it is
shown that all Osserman algebraic curvature tensors in $\Ro$ have a Clifford
structure, and Proposition~\ref{p2} which says that a Riemannian manifold with
such a curvature tensor is flat or rank-one symmetric. Proposition~\ref{p1} is
proved in Section~\ref{sec:3} modulo two Lemmas (their proofs are given in
Section~\ref{s:lemmas}). Proposition~\ref{p2} is proved in Section~\ref{s:cliff8}.

\section{Manifolds with Clifford structure. Proof of the Theorem}
\label{sec:2}

We follow the two-step approach to the Osserman Conjecture
suggested in \cite{GSV}:
\begin{itemize}
\item[(i)] find all Osserman algebraic curvature tensors;
\item[(ii)] classify Riemannian manifolds having curvature tensor as in
\textup{(i)}.
\end{itemize}

The standard tool for \textup{(ii)} is the second Bianchi identity,
although the proof can sometimes be quite technically involved.

The difficult part is \textup{(i)}, but thanks to the remarkable
construction of \cite{G,GSV}, we know the right candidate for \textup{(i)},
a typical Osserman algebraic curvature tensor:

\begin{definition}\label{d3}
An algebraic curvature tensor $R$ in $\Rn$ has a {\it $\Cliff(\nu)$-structure}
($\nu \ge 0$), if there exist anticommuting skew-symmetric orthogonal operators
$J_1, \dots , J_\nu$, and the numbers $\la_0, \mu_1, \ldots \mu_\nu$,
with $\mu_s \ne \la_0$, such that
\begin{multline}\label{eq:Clstr}
R(X, Y) Z = \la_0 (\< X, Z \> Y - \< Y, Z \> X) \\
+ \sum\nolimits_{s=1}^\nu \tfrac 13 (\mu_s - \la_0) (2 \< J_sX, Y \> J_sZ +
\< J_sZ, Y \> J_sX - \< J_sZ, X \> J_sY).
\end{multline}
A Riemannian manifold $M^n$ has a {\it $\Cliff(\nu)$-structure} if its
curvature tensor at every point does.
\end{definition}

The fact that skew-symmetric operators $J_s$ are orthogonal and
anticommute is equivalent to each of the following sets of equations:
$\<J_sX, J_qX\> = \delta_{sq}\|X\|^2$ and $J_sJ_q+J_qJ_s =-2 \delta_{sq} I_n$,
for all $s,q =1, \ldots, \nu$ and all $X \in \mathbb R^n$.

The Jacobi operator of the algebraic curvature tensor $R$ with the Clifford
structure given by \eqref{eq:Clstr} has the form
\begin{equation} \label{e2}
R_XY =  \la_0 (\| X\|^2 Y - \< Y, X \> X) +
\sum\nolimits_{s=1}^\nu (\mu_s -\la_0) \< J_sX, Y \> J_sX,
\end{equation}
and the tensor $R$ can be reconstructed from \eqref{e2} using polarization
and the first Bianchi identity.

It follows that $\Cliff(\nu)$ algebraic curvature tensor (manifold) is Osserman
(pointwise Osserman, respectively). Indeed, for any unit vector $X$, the Jacobi
operator $R_X$ given by \eqref{e2}
has constant eigenvalues $\la_0, \la_1, \dots , \la_{k-1}$, where $\la_1, \dots , \la_{k-1}$
are the $\mu_s$'s without repetitions. The eigenspace corresponding to
the eigenvalue $\la_\a, \; \a \ne 0$, is
$E_{\la_\a}(X) = \Span_{s: \mu_s = \la_\a} (J_sX)$, and the $\la_0$-eigenspace
is $E_{\la_0}(X) = (\Span (X, J_1X, \dots , J_\nu X))^\perp$, provided
$\nu < n-1$.

Following the above approach, we deduce the Theorem
from two Propositions:

\begin{proposition}
\label{p1}
An Osserman algebraic curvature tensor in $\Ro$ has a Clifford structure.
\end{proposition}

\begin{proposition} \label{p2}
A pointwise Osserman manifold with a
Clifford structure is flat or rank-one symmetric.
\end{proposition}

\section{Proof of Proposition~\ref{p1}}
\label{sec:3}

Let $R$ be an Osserman algebraic curvature tensor in $\Ro$, with
the Jacobi operator having $k$ distinct eigenvalues
$\la_0, \la_1, \ldots, \la_{k-1}$ whose multiplicities
are $m_0, m_1, \dots, m_{k-1}$, respectively, with $m_0+m_1 + \dots \linebreak[2]+ m_{k-1} = 7$.
For a nonzero vector $X$, the eigenvalues of $R_X$ are then
$0, \la_0 \|X\|^2, \la_1 \|X\|^2, \ldots, \la_{k-1} \|X\|^2$, and the
corresponding eigenspaces are $\Span(X)$ and
$E_{\la_\a \|X\|^2}(X)= E_{\la_\a}(X/|X\|)$, respectively.
(Recall that ``the eigenvalues of $R_X$" are just the eigenvalues of
the operator $R_X$ acting in $\Ro$, and there are $8$ of them, counting the multiplicities.
When we say ``the eigenvalues of the Jacobi operator", we do not count the eigenvalue $0$
corresponding to the vector $X$, so the sum of multiplicities is $7$. We allow one of the
$\lambda_\alpha$'s to be zero, in which case $\Ker R_X = \Span(X) \oplus E_0(X)$).

If $k=1$, there is nothing to prove: the curvature is constant. If $k=2$, the claim follows from
{\cite[Proposition 1]{Ndga}}:
\begin{lemma} \label{l:2eigen}
An Osserman algebraic curvature tensor in $\Rn, \; n \ne 16$, with
the Jacobi operator having two distinct eigenvalues, has a Clifford structure
$\Cliff(\nu)$, with $\nu < n/2$.
\end{lemma}

We may assume therefore, that $k \ge 3$. 
Let $m_0$ be the maximal multiplicity, and $\nu = 7 - m_0$ be the sum
of all the others. Two cases are possible:
\begin{enumerate}
  \item [(A)] either each of the $m_1, m_2, \ldots, m_{k-1}$ is $1$ or $2$,
  \item [(B)] or $k=3$ and the multiplicities are $m_0=3, m_1=1, m_2=3$, up to
relabelling.
\end{enumerate}
To prove that the Osserman algebraic curvature tensor $R$ has a Clifford structure,
it is sufficient to show that for every eigenvalue $\la_\a \; (\a > 0)$ of
the Jacobi operator, there exist $m_\a$ anticommuting
skew-symmetric orthogonal operators $J^\a_1, \dots , J^\a_{m_\a}$ such that
for every unit vector $X$, the eigenspace $E_{\la_\a}(X)$ of $R_X$ is spanned
by $J^\a_1X, \dots , J^\a_{m_\a}X$. Indeed, if this is the case, then
for every unit vector $X$ and every $Y \in \Rn$,
\begin{equation*}
R_XY =  \la_0 (\| X\|^2 Y - \< Y, X \> X) +
\sum\nolimits_{\a=1}^p \sum\nolimits_{s=1}^{m_\a} (\la_\a -\la_0) \< J^\a_sX, Y \> J^\a_sX,
\end{equation*}
as the symmetric operators on the left-hand side and on the right-hand side have the same
eigenvalues and eigenspaces. So the Jacobi operator has the form \eqref{e2}, up to
relabelling through: $J_1 = J^1_1, J_2=J^1_2, \linebreak[2]\ldots , J_\nu = J^p_{m_p}$. The only remaining thing to
check is that the operators $J^\a_s,\, J^\beta_q$, with $\beta \ne \alpha$ anticommute,
which easily follows from the fact that for all unit vectors
$X,\quad J^\a_sX \in E_{\la_\a}(X) \perp E_{\la_\beta}(X) \ni J^\beta_qX$.

In the both cases (A) and (B) above, the proof will follow from the two Lemmas below.

First of all, as is immediate from Definition~\ref{d3}, shifting an algebraic curvature tensor
with a Clifford structure by a constant curvature tensor results in an algebraic
curvature tensor, still having a Clifford structure. So we can always assume that
$\la_0$, the eigenvalue with the highest multiplicity, is zero.

Next, for all the eigenvalues of multiplicity $1$ or $2$, we can find polynomial
vectors spanning the corresponding eigenspaces, with the properties similar to
those of the $J_sX$'s:

\begin{lemma} \label{l:mult1,2}
\emph{1.} Let $R$ be an Osserman algebraic curvature tensor in $\Rn$ whose Jacobi
operator has $k$ distinct eigenvalues, one of which is zero, and let $\la \ne 0$ be a
simple eigenvalue. Then for every $X \ne 0$, the eigenspace $E_{\lambda\|X\|^2}(X)$ of $R_X$ is
spanned by a vector $P(X)$ all of whose components are odd homogeneous polynomial
of degree $2m+1 \le k-1$, and for all unit vectors $X \in \Rn$,
$$
\< P(X), X \> = 0, \qquad \| P(X) \|^2 = 1, \qquad P(P(X))= - X.
$$

\emph{2.} Let $R$ be an Osserman algebraic curvature tensor in $\Rn$ whose Jacobi
operator has $k$ distinct eigenvalues, one of which is zero, and let $\la \ne 0$
be an eigenvalue of multiplicity $2$. Then for every $X \ne 0$, the eigenspace
$E_{\lambda\|X\|^2}(X)$ of $R_X$ is spanned by vectors $U(X), V(X)$ all of whose
components are odd homogeneous polynomial of degree $2m+1 \le k-1$, and
for all unit vectors $X \in \Rn$,
$$
\< U(X), X \> = \< U(X), X \> = \< U(X), V(X) \> = 0,
\qquad \| U(X) \|^2 = \| V(X) \|^2 = 1.
$$
\end{lemma}

Note that in the case when the degree of $P$ (respectively, $U,V$) is $1$, the
equations of Lemma \ref{l:mult1,2} simply mean that the corresponding linear operators
in $\Rn$ are orthogonal, skew-symmetric and anticommute.

The reason we require the multiplicity to be less than or equal to two is topological:
any $2$-dimensional plane bundle over the real projective space $\mathbb RP^n$ of
dimension $n\ge 3$ is a sum of line bundles \cite{Lev,H}, which is no longer true
for $3$-dimensional bundles.

\vskip .3cm
\textit{Case A}. The Jacobi operator has three eigenvalues $0, \lambda_1, \lambda_2$
of multiplicities $3, 1, 3$, respectively. From the first assertion of Lemma~\ref{l:mult1,2},
the eigenspace $E_{\la_1 \|X\|^2}(X)$ is spanned by a homogeneous polynomial vector $P(X)$
of degree $1$, that is, there exists an orthogonal skew-symmetric operator $J$ such
that $P(X) = JX$ and $R_XJX = \la_1 \|X\|^2JX$, for all $X \in \Ro$. Introduce an algebraic
curvature tensor $\hat R$ by
\begin{equation*}
\hat R(X, Y) Z = R(X, Y) Z -
\tfrac {\la_1}3 (2 \< JX, Y \> JZ + \< J Z, Y \> J X - \< J Z, X \> J Y).
\end{equation*}
Its Jacobi operator has the form $\hat R_X Y = R_X Y - \la_1 \< JX, Y \> JX$. For
every unit vector $X, \; \hat R_XJX =0$, and the restrictions of $\hat R_X$ and
$R_X$ to $(JX)^\perp$ coincide, so the Jacobi operator $\hat R_X$ has two eigenvalues,
$\la_2$ and $0$, with the eigenspaces $\hat E_{\la_2}(X)=E_{\la_2}(X)$,
of dimension $3$, and $\hat E_0(X) = \Span(JX)\oplus E_{0}(X)$, of dimension $4$,
respectively. It follows that the algebraic
curvature tensor $\hat R$ is Osserman, with the Jacobi operator having two eigenvalues.
By Lemma~\ref{l:2eigen}, $\hat R$ has a Clifford structure which is a
$\Cliff(3)$-structure, and so
there exist anticommuting orthogonal skew-symmetric operators $J_1, J_2, J_3$ in $\Ro$
such that for any unit vector $X,\; \hat E_{\la_2}(X)=\Span(J_1X, J_2X, J_3X)$.

Hence $R$ has a Clifford structure, as for any unit vector
$X,\; E_{\la_1}(X)=\Span(JX)$ and $E_{\la_2}(X)= \hat E_{\la_2}(X)=\Span(J_1X, J_2X, J_3X)$.

\vskip .3cm
\textit{Case B}. Let all the multiplicities $m_1, \ldots , m_{k-1}$ be $1$ or $2$.
From Lemma~\ref{l:mult1,2} it follows that for every $X\ne 0$, all the eigenspaces of the
operator $R_X$ corresponding to the nonzero eigenvalues are spanned by polynomial vectors.
Label them through $P_1(X), \ldots, P_\nu(X), \; \nu =m_1+ \ldots + m_{k-1}$. The proof of the
Proposition is then completed with the following Lemma (the fact that the $P_s(X)$'s verify
the assumptions easily follows from Lemma~\ref{l:mult1,2}).

\begin{lemma} \label{l:polys to Cl}
Let $R$ be an Osserman algebraic curvature tensor in $\Ro$
whose Jacobi operator has $\nu \le 6$ nonzero eigenvalues
$\mu_1, \mu_2, \ldots, \mu_\nu$, counted with multiplicities. Assume that the
multiplicity of each of the $\mu_s$'s is not greater than $\min\{2, 7 - \nu\}$.
Let $P_1(X), P_2(X), \ldots, P_\nu(X)$ be eigenvectors of
$R_X$ such that for every $s = 1, \ldots , \nu$, all the components of
$P_s(X)$ are homogeneous polynomials of odd degree $d_s$, and
\begin{gather}
R_X P_s(X) = \mu_s \|X\|^2 P_s(X), \label{eq:P_i1}
\\
\<P_s(X), P_t(X)\> = \delta_{st} \|X\|^{2d_s}, \qquad \<P_s(X), X\> = 0,
\label{eq:P_i2}
\end{gather}
for all $X \in \Ro$ and all $s,t = 1, \ldots , \nu$, and additionally, in the case $\nu =
6$,
\begin{equation}\label{eq:pp}
  P_s(P_s(X)) = - \|X\|^{d_s^2 -1} X
\end{equation}
for all $X \in \Ro$ and all $s = 1, \ldots , 6$.

Then $R$ has a $\Cliff(\nu)$-structure.
\end{lemma}

\section{Proof of Lemma~\ref{l:mult1,2} and Lemma~\ref{l:polys to Cl}}
\label{s:lemmas}

In this Section, we give the proofs of Lemma~\ref{l:mult1,2} (Subsection~\ref{ss:mult1,2})and
Lemma~\ref{l:polys to Cl} (Subsection~\ref{ss:polys to C1}). The latter one will use some
algebraic results from Subsection~\ref{ss:ca}.
Throughout this Section, for $\mathbb{K} =\mathbb R$ or $\mathbb C, \; \mathbb{K}[X]$ is the
ring of polynomials in $X =(x_1, \ldots, x_n)$ (or in $X =(x_1, \ldots, x_8)$,
depending on the context) over $\mathbb{K}$. The polynomial $\|X\|^2 = \sum_j x_j^2$
is irreducible in $\mathbb{K}[X]$ (for $n \ge 3$), and $\mathbb{K}[X]$ is a unique factorization
domains (UFD). We will say that a matrix (a vector) with entries from $\mathbb{K}[X]$
is divisible by $\|X\|^2$, if each of its entries is.

\subsection{Proof of Lemma~\ref{l:mult1,2}}\label{ss:mult1,2}

The first assertion of Lemma \ref{l:mult1,2} follows from Lemma 2.1 of \cite{Nhjm}.
Here we prove the second one. Throughout the proof, we will freely switch between
operators and their matrices using the same notation. We call a finite set of polynomials
coprime, if their greatest common divisor is one.

Let $\la_0 =0, \la_1 = \la, \la_2, \dots, \la_{k-1}$ be the eigenvalues of the Jacobi
operator.

For every $X \in \Rn$, define a symmetric operator $W_X: \Rn \to \Rn$ by
$W_X = \rho^{-1} R_X \prod_{\a =2}^{k-1} (R_X - \la_\a \|X\|^2 I_n)$,
where
$\rho = \la \prod_{\a=2}^{k-1} (\la - \la_\a)$. The operator $W_X$ is an
orthogonal projection on the $2$-plane $E_{\la \|X\|^2}(X)$ followed by multiplication
by $\|X\|^{2k-2}$, and the entries of its matrix are homogeneous polynomials
of even degree $2k-2$ of the coordinates $x_1, \dots , x_n$ of vector $X$, so that
\begin{equation*}
W_X^2 = \|X\|^{2k-2} W_X, \quad \text{and} \quad \rk \,W_X = 2 \quad \text{when}
\quad X \ne 0.
\end{equation*}
Let $q(X)$ be the greatest common divisor of the entries of the matrix
$W_X$ over $\mathbb R[X]$. Then $q(X)$ divides $\Tr W_X = 2 \|X\|^{2k-2}$,
and so $q(X) = \|X\|^{2k-2-2s}$ for some $s \le k-1$.
Dividing $W_X$ by $q(X)$ we get a symmetric polynomial matrix $S_X$ whose entries
are coprime homogeneous polynomials of degree $2s$ satisfying
\begin{equation}\label{eq:s rk2}
S_X^2 = \|X\|^{2s} S_X, \qquad \qquad \rk \,S_X = 2 \quad \text{when}
\quad X \ne 0,
\end{equation}
and the two-dimensional $\|X\|^{2s}$-eigenspace of $S_X$
is $E_{\la \|X\|^2}(X)$.

Since for any unit vector $X,\; E_\la(X) = E_\la(-X)$, the two-plane bundle $E_\la$
over the unit sphere in $\Rn$ descends to a two-plane bundle over $\mathbb RP^{n-1}$.
Any such bundle is a sum of two line bundles when $n \ge 3$ \cite{Lev,H},
so we can choose two orthonormal vector fields $s_1(X), s_2(X)$ spanning
$E_\la(X)$ for $X \in S^{n-1}$ which are odd with respect to the antipodal
map.

Let $\tilde s_1, \tilde s_2$ be vector fields in $\Rn$ whose components are
odd homogeneous polynomials of the same degree and
whose restrictions to $S^{n-1}$ are uniformly close to $s_1, s_2$, respectively.
As the polynomial matrix $S_X$ is a projection on the $2$-plane $E_\la(X/\|X\|)$
followed by multiplication by $\|X\|^{2s}$, the components of vector fields
$$
P_1(X) = S_X \tilde s_1(X), \quad P_2(X) = S_X \tilde s_2(X)
$$
are odd homogeneous polynomials of the same degree, with $P_1(X)$ and $P_2(X)$
pointwise linearly independent on the sphere $S^{n-1}$ (and hence in $\Rn \setminus 0$),
and spanning $E_{\la \|X\|^2}(X)$.
Replacing $P_1$ by $\|P_2\|^2 P_1 - \<P_1, P_2\> P_2$, we can assume that
for all $X, \; P_1(X) \perp P_2(X)$. Dividing each of
$P_1(X), P_2(X)$ by the greatest common divisor of its components,
we get two homogeneous polynomial vectors $P_1(X), P_2(X)$,
each with coprime components, and such that
\begin{equation}\label{eq:sxpi}
\<P_1(X), X\> = \<P_2(X), X\> = \<P_1(X), P_2(X)\> = 0, \quad
S_X P_i(X) = \|X\|^{2s} P_i(X),
\end{equation}
for all $X \in \Rn,\; i = 1,2$. The polynomials $a(X) = \|P_1(X)\|^2$
and $b(X) = \|P_2(X)\|^2$ are nonzero outside the origin in $\mathbb R^n$,
and both $P_1$ and $P_2$ are still odd, as from $P_1 (X), P_2(X) \perp X$
it follows that the self-maps of $S^{n-1}$ defined by
$X \to a(X)^{-1/2}P_1(X), \; X \to b(X)^{-1/2}P_2(X)$
are homotopic to the identity.

From \eqref{eq:s rk2} and \eqref{eq:sxpi}, we get the following matrix equation:
$$
(P_1 | P_2)
\left(
\begin{matrix} a^{-1} & 0 \\ 0 & b^{-1} \end{matrix}
\right)
(P_1 | P_2)^t =  \|X\|^{-2s} S_X,
$$
for all $X \ne 0$.
Let $a = a'c, \, b = b'c$, with $a',\, b'$ coprime polynomials. Then
$$
\|X\|^{2s} b'(X) P_1(X) P_1(X)^t \!\!+ \|X\|^{2s} a'(X) P_2(X) P_2(X)^t \!\!=
a'(X) b'(X) c(X) S_X,
$$
so the polynomial matrix
$\|X\|^{2s} b'(X) P_1(X) P_1(X)^t$ is divisible by $a'(X)$. Since $a'$ and
$b'$ are coprime and the components of $P_1(X)$ are coprime, $\|X\|^{2s}$
must be divisible by $a'(X)$, so that $a'(X) =
\|X\|^{2m}$, for some $m \le s$. Similarly, $b'(X) =\|X\|^{2m'}$. Again,
since $a'$ and $b'$ are coprime, no more than one of the numbers $m$,
$m'$ can be nonzero. Let say $m'=0$, that is, $b'(X) = 1$. Then
$$
\|X\|^{2s-2m} P_1(X) P_1(X)^t + \|X\|^{2s} P_2(X) P_2(X)^t = c(X) S_X.
$$
Since the polynomial $\|X\|^2$ is irreducible, and the entries of $S_X$ are
coprime, $c(X) = \|X\|^{2s-2m} d(X)$ for some polynomial $d$, and we come to
$$
P_1(X) P_1(X)^t + \|X\|^{2m} P_2(X) P_2(X)^t = d(X) S_X.
$$
The number $m$ must be even. Indeed, both $P_1(X)$ and $P_2(X)$
are  of odd degree, so the polynomials $a(X) = \|P_1(X)\|^2$ and
$b(X) = \|P_2(X)\|^2$ are of degree $2 \mod 4$. But $a = a' c = \|X\|^{2s}
d(X), \; b = c = \|X\|^{2s-2m} d(X)$, so $2m$ is divisible by $4$.

Then $\|X\|^m$ is a polynomial and we can introduce polynomial vectors
$P = P_1, \,\, Q = \|X\|^m P_2$ which satisfy
$P(X) P(X)^t + Q(X) Q(X)^t = d(X) S_X$, and so for every $Y \in \mathbb R^n$,
\begin{equation}\label{eq:24}
\<P(X), Y\>^2 + \<Q(X), Y\>^2 = d(X) \<S_XY, Y\>.
\end{equation}

The components of $P(X)$ are still coprime over $\mathbb R [X]$,
while the greatest common divisor of the components of $Q(X)$ is
$\|X\|^m$. We want to find two polynomial vectors $U(X), V(X)$ spanning the
same subspace as $P(X), Q(X)$ (which is the eigenspace $E_{\la\|X\|^2}(X)$ of $R_X$)
and satisfying \eqref{eq:24} with $d(X) = 1$.

Let $r(Z)$ be an irreducible factor of $d(Z)$ over
$\mathbb C[Z]$. Then either $r(Z)$ is essentially real, that is, $\bar r (Z) = r(Z)$
up to multiplication by a nonzero complex number, or $d(Z)$ is also divisible by $\bar r(Z)$
(where $\bar f$ is defined by $\bar f(Z) = \overline{f(\bar Z)}$).
In both cases, \eqref{eq:24} implies that for every $Y \in \mathbb R^n$,
the complex polynomial $\<P(Z) + i Q(Z), Y\> \<P(Z) - i Q(Z), Y\>$ is divisible
by $r(Z)$ over $\mathbb C [Z]$ (the scalar product is Euclidean, not Hermitian). Let
say $\<P(Z) + i Q(Z), Y\>$ be divisible by $r(Z)$, for all real $Y$, so that
$r(Z)$ divides every component of the vector $P(Z) + i Q(Z)$.

In the first case, since $\bar P = P, \, \bar Q = Q, \, \bar r = r$, the polynomial $r(Z)$
also divides every component of $P(Z) - i Q(Z)$. Then every component of $P(Z)$
is divisible by $r(Z)$ over $\mathbb C[Z]$, and so every component of $P(X)$ is
divisible by $r(X)$ over $\mathbb R[X]$. This contradicts to the fact that the
components of $P$ are coprime over $\mathbb R[X]$.

In the second case, let $F(Z)$ be a polynomial vector such
that $P(Z) + i Q(Z) = r(Z) F(Z)$. For real $Z = X$, let
$r(X) = r_1(X) + i r_2(X), \; F(X) = A(X) + i B(X)$, with $r_1, r_2$ real
polynomials, and $A, B$ real polynomial vectors. Then
$A = (P r_1 + Q r_2) (r_1^2 + r_2^2)^{-1}, \;
B = (-P r_2 + Q r_1) (r_1^2 + r_2^2)^{-1}$, and so
\begin{multline*}
A(X) A(X)^t + B(X) B(X)^t
= (r_1(X)^2 + r_2(X)^2)^{-1} (P(X) P(X)^t + Q(X) Q(X)^t)  \\
= ((r_1(X)^2 + r_2(X)^2)^{-1} d(X)) S_X.
\end{multline*}
Since $d(Z)$ is divisible by irreducible polynomials $r(Z)$ and $\bar r(Z)$ over
$\mathbb C[Z]$, it is also divisible by $r(Z) \bar r(Z)$.
Then, over $\mathbb R[X]$, the real polynomial $d(X)$ is divisible by the real
polynomial $r(X) \bar r(X) = r_1(X)^2 + r_2(X)^2$. Hence
$(r_1(X)^2 + r_2(X)^2)^{-1} d(X)$ is a polynomial, of degree less than $\deg d(X)$.

Repeating this procedure, we finally come to two polynomial vectors,
$U(X)$ and $V(X)$, such that
\begin{equation}\label{eq:25}
 U(X) U(X)^t + V(X) V(X)^t = S_X.
\end{equation}
It remains to check that $U(X), V(X)$ is the sought pair of polynomial vectors.

Since by \eqref{eq:s rk2} $\rk S_X = 2$ for all $X \ne 0$,
the vectors $U(X), V(X)$ must be pointwise linearly independent outside the origin.
Substituting $S_X$ from \eqref{eq:25} to \eqref{eq:s rk2} we get
$U(X) ((\|U(X)\|^2 - \|X\|^{2s})U(X)^t +\<U(X),V(X)\>V(X)^t)+ V(X) ((\|V(X)\|^2 - \|X\|^{2s})V(X)^t +\<U(X),V(X)\>U(X)^t) = 0$
which implies, by the linear independence of $U(X)$ and $V(X)$, that
$$
\|U(X)\|^2 = \|V(X)\|^2 = \|X\|^{2s}, \quad \<U(X), V(X) \> = 0.
$$
It also follows from \eqref{eq:25}, that all the components of $U(X)$ and
$V(X)$ are homogeneous polynomials of degree $s$, and for $X \ne 0$,
the vectors $U(X), V(X)$ span $E_{\la\|X\|^2}(X)$. In
particular, as $E_{\la\|X\|^2}(X) \perp X$,
$$
\<U(X), X \> = \<V(X),  X \> = 0,
$$
and so $s = \deg U(X) = \deg V(X)$ is odd (otherwise the homotopy $\cos \T X +
\sin \T U(X)$ joins the identity map of $S^{n-1}$ with the one of even
topological degree). Finally, by construction of $S_X, \; s \le k-1$.

\subsection{Two algebraic facts}\label{ss:ca}

Let for $X = (x_1, \ldots, x_8), \; (\|X\|^2)$ be the ideal of $\mathbb R[X]$
generated by $\|X\|^2 = \sum_j x_j^2$, and $\mathbf{R} = \mathbb R[X]/(\|X\|^2)$, with
$\pi:\mathbb R[X] \to \mathbf{R}$ the natural projection.
The field of fractions $\mathbb{F}$ of the ring $\mathbf{R}$ is isomorphic to
the field $\mathbb{L}_7 = \mathbb{R}(x_1, \ldots, x_7, \sqrt{-d})$, where
$d=x_1^2 + \ldots + x_7^2$
(an isomorphism from $\mathbb{L}_7$ to $\mathbb{F}$ is induced by the map
$(a + b \sqrt{-d})/c \to (a + bx_8)/c$, with
$a, b, c \in \mathbb{R}[x_1, \ldots, x_7], \; c \ne 0$).

The \textit{level} $s(\Bbbk)$ of a field $\Bbbk$ is the minimal number of elements
$c_1, \ldots , c_m \in \Bbbk$ such that $c_1^2 + c_2^2 + \ldots + c_m^2 = -1$ \cite{Pf}.
If $-1$ is not a sum of squares in $\Bbbk$, then $s(\Bbbk) = \infty$. By Pfister's
Theorem, the level of any field is a power of $2$, if finite.

We need the following two facts:
\begin{fact}[\cite{Nag}]
\label{fact:1}
The rings $\mathbf{R}$ and $\mathbf{R}\otimes{\mathbb{C}} = \mathbb C[X]/(\|X\|^2)$ are unique
factorization domains \emph{(}UFD's\emph{)}.
\end{fact}
\begin{fact}[{\cite[Theorem 3.1.4]{Pf}}]\label{fact:2}
The level of the field $\mathbb{F} \, (\,\cong \mathbb{L}_7)$ is $4$.
\end{fact}

\begin{remark}
Fact~\ref{fact:2} is equivalent to the following: if a sum of squares
of no more than four elements in $\mathbf{R}$ is zero, than all the elements are
zeros. Indeed, if $a_1^2+\ldots+a_m^2 =0$ for $a_1, \ldots , a_m \in \mathbf{R}$,
with $m \le 4$, and say $a_m \ne 0$, then in
$\mathbb{F}, \quad (a_1/a_m)^2 + \ldots + (a_{m-1}/a_m)^2 = -1$
which is a contradiction, as $m-1 < 4$.

Equivalently, if a sum of squares of no more than four real polynomials in
eight variables $x_1, \ldots, x_8$ is divisible by $\|X\|^2$, then each of them
is. In the proof of Proposition~\ref{p2} (Lemma~\ref{l:octonions}), we will
also use a similar fact for polynomials in six variables: if a sum of squares of
two real polynomials in $x_1, \ldots , x_6$ is divisible by $\sum_{i=1}^6 x_i^2$,
then both are divisible by $\sum_{i=1}^6 x_i^2$ (this also follows
from \cite[Theorem 3.1.4]{Pf}: the level of the field
$\mathbb{L}_5 = \mathbb{R}(x_1, \ldots, x_5, \sqrt{-(x_1^2 + \ldots + x_5^2)}\,)$
is $4$).
\end{remark}

In the proof of Lemma~\ref{l:polys to Cl}, we will use the following Lemma for
polynomial matrices.
\begin{lemma} \label{l:4squares}
Let $A$ be an $8 \times m$ polynomial matrix \emph{(}with entries in $\mathbb{R}[X]$\emph{)},
and $D$ be a real diagonal $m \times m$ matrix with positive diagonal entries.
Assume that $ADA^t$ is divisible by $\|X\|^2$.
Then:
\begin{enumerate}
  \item[\emph{(i)}] If $m \le 4$, then $A$ itself is divisible by $\|X\|^2$.
  \item[\emph{(ii)}]  If $m = 5$, then there exist polynomial vectors $U$ and $Q$ of dimension
$8$ and $5$, respectively, and an $8 \times 5$ polynomial matrix $B$ such that
\begin{equation}\label{e:rk1+x^2}
  A=UQ^t + \|X\|^2 B, \quad \text{and the polynomial $Q^tDQ=\<DQ,Q\>$ is divisible by $\|X\|^2$}.
\end{equation}
\end{enumerate}
\end{lemma}

\begin{proof}
Denote $d_1, \ldots, d_m >0$ the diagonal entries of the matrix $D$. Let $a$
be an $8 \times m$ matrix with entries in $\mathbf{R}$, such that
$a_{js} = \pi(\sqrt{d_s}A_{js}), \; s = 1, \ldots, m,\;j = 1, \ldots , 8$.
Then by assumption, $aa^t = 0$ which implies that for any two rows
$a_i, a_j$ of the matrix $a,\quad \|a_i\|^2 = \|a_j\|^2 = \< a_i, a_j\> = 0$
(where $\<\cdot,\cdot\>$ and $\|\,\cdot\,\|^2$ in the free module $\mathbf{R}^m$ have an
obvious meaning).
\vskip .2cm
(i) If $m \le 4$, then, for every $j, \; \|a_j\|^2$ is a sum of no more than
four squares in $\mathbf{R}$, and so the equation $\|a_j\|^2 = 0$ implies
$a_j = 0$ by Fact~\ref{fact:2}.
\vskip .2cm
(ii) Let $m=5$. We have to find $u \in \mathbf{R}^8,\,q \in \mathbf{R}^5$ such
that $a = uq^t$ (then automatically $\|q\|^2=0$, unless $a=0$, in which case
there is nothing to prove).

Without loss of generality, assume that $a_{11}\ne 0$. Denote $u_1$ the
greatest common divisor of the entries in the first row of $a$
($\mathbf{R}$ is a UFD by Fact~\ref{fact:1}), and let $a_{1s} =u_1q_s, \; s = 1, \ldots, 5$,
with the $q_s$'s coprime (note that none of the $q_s$'s is zero, because otherwise
$0 = \|a_1\|^2$ is a sum of no more than four squares, and so $a_1 =0$).

As for any two rows $a_i, a_j \quad \|a_i\|^2 = \|a_j\|^2 = \< a_i, a_j\> = 0$,
it follows that for all $j =2, \ldots, 8, \quad \|a_{11} a_j-a_{j1}a_1\|^2 = 0$
which is a sum of no more than four squares in $\mathbf{R}$.
So by Fact~\ref{fact:2}, $a_{11} a_j=a_{j1}a_1$, that is,
$q_1 a_{js} = a_{j1} q_s$ for all $j \ge 2, s$. As $\mathbf{R}$ is a UFD, there
exist $u_2, \ldots, u_8$ such that $a_{js} = u_jq_s$, for all $j, s$.
\end{proof}

\subsection{Proof of Lemma~\ref{l:polys to Cl}}\label{ss:polys to C1}

Without loss of generality, we can assume that none of the
polynomial vectors $P_s(X)$ is divisible by $\|X\|^2$ (for every $s$, at least
one of the components of $P_s(X)$ is not divisible by $\|X\|^2$).

Let $d = \max\{d_1, \ldots, d_\nu\}$. Multiplying the vectors $P_s(X)$,
with $d_s < d$, by an appropriate power of $\|X\|^2$, we can assume that all the components
of all the $P_s$'s have degree $d$, and at least one of the $P_s(X)$'s is still not
divisible by $\|X\|^2$.

We want to show that $d =1$, that is, each $P_s$ is a linear operator in $\Ro$.
Assume that $d >1$ (recall that all the $d_s$'s are odd, and so is $d$).

Let $A(X)$ be $8 \times \nu$ polynomial matrix whose
columns are $P_1(X), \ldots, P_\nu(X)$, and let $\Lambda$ be $\nu \times \nu$
constant diagonal matrix with diagonal entries $\mu_1, \ldots, \mu_\nu$. Then
by \eqref{eq:P_i1} the symmetric $8 \times 8$ matrices $A(X) \Lambda A(X)^t$ and
$\|X\|^{2d-2} R_X$ have the same eigenspaces and eigenvalues for all $X \ne 0$, so
\begin{equation}\label{eq:plp=rx}
A(X) \Lambda A(X)^t = \|X\|^{2d-2} R_X, \qquad A(X)^t A(X) = \|X\|^{2d} I_8,
\end{equation}
(the second equation follows from \eqref{eq:P_i2}). Raising both sides of the first
equation of
\eqref{eq:plp=rx} to the power $k < d$ and using the second one we get
\begin{equation}\label{eq:plkp=rxk}
A(X) \Lambda^k A(X)^t = \|X\|^{2(d-k)} R_X^k, \qquad k = 1, \ldots , d-1,
\end{equation}
so all the polynomial matrices $A(X) \Lambda^k A(X)^t$ with $k \le d-1$ are divisible
by $\|X\|^2$. We are going to prove that the matrix $A(X)$ is itself divisible by
$\|X\|^2$. This will lead to a contradiction with the fact that at least one of the
$P_s(X)$'s is not divisible by $\|X\|^2$.

Let $\nu \le 4$. By \eqref{eq:plkp=rxk}, the matrix $A(X) \Lambda^2 A(X)^t$ is
divisible by $\|X\|^2$, with $\Lambda^2$ diagonal, positive definite matrix. Then
by Lemma \ref{l:4squares}(i), $A(X)$ is divisible by $\|X\|^2$.

Let $\nu = 5$. Without loss of generality, assume that the $\mu_s$'s are labelled in
such a way that $\mu_1^{-1} \le \mu_2^{-1} \le \ldots \le \mu_5^{-1}$.
Let $m (=1$ or $2)$ be the multiplicity of $\mu_1$, and let $A(X) = \bigl(A_1 \,|\, A_2 \bigr)$,
where $A_1, A_2$ are submatrices of $A(X)$ consisting of the first $m$, and the last
$5-m$ columns, respectively. Denote $\widetilde{\Lambda}$ the $(5-m)\times(5-m)$
diagonal matrix whose diagonal entries are $\mu_{m+1} - \mu_1, \ldots, \mu_5-\mu_1$.
Then from \eqref{eq:plp=rx} and \eqref{eq:plkp=rxk}, with $k =2$, both matrices
\begin{equation}\label{eq:tildelambda}
\mu_1 A_1A_1^t + A_2\widetilde{\Lambda} A_2^t, \quad\text{and}\quad
\mu_1^2 A_1A_1^t + A_2\widetilde{\Lambda}^2 A_2^t
\quad \text{are divisible by $\|X\|^2$},
\end{equation}
and then, so is $A_2(\widetilde{\Lambda}^2 - \mu_1\widetilde{\Lambda}) A_2^t$.
By Lemma \ref{l:4squares}(i), this implies that $A_2$ is
divisible by $\|X\|^2$, as $5-m \le 4$ and the diagonal matrix
$\widetilde{\Lambda}^2 - \mu_1\widetilde{\Lambda}$ is definite: all its diagonal
entries $\mu_s^2 - \mu_1\mu_s = \mu_s^2 \mu_1 (\mu_1^{-1} - \mu_s^{-1})$ have the
same sign. Then $\mu_1A_1A_1^t$ is divisible by $\|X\|^2$ by \eqref{eq:tildelambda},
and so is $A_1$, again by Lemma \ref{l:4squares}(i).

Let now $\nu = 6$. In this case, by the assumption, all the $\mu_s$'s are distinct.
We can assume that they are labelled in
such a way that $\mu_1^{-1} < \mu_2^{-1} < \ldots < \mu_6^{-1}$.

First consider the case $d \ge 5$. Then by \eqref{eq:plkp=rxk}, all four polynomial
matrices $A(X) \Lambda^k A(X)^t, \; k =1, 2, 3, 4$ are divisible by $\|X\|^2$. Set
$f(t) = t(t-\mu_1)(t-\mu_2)(t-\mu_3)$, and denote
$\widetilde{\Lambda}= \mathrm{diag} \{f(\mu_4), f(\mu_5), f(\mu_6)\}$. Note that
all three diagonal entries $f(\mu_4), f(\mu_5), f(\mu_6)$ of $\widetilde{\Lambda}$
have the same sign. Let $A(X) = \bigl(A_1 \,|\, A_2 \bigr)$, where $A_1$ and $A_2$ consist
of the first three and the last three columns of $A(X)$, respectively.
Then the matrix $A(X) f(\Lambda) A(X)^t = A_2 \widetilde{\Lambda} A_2^t$
is divisible by $\|X\|^2$, as $f(\Lambda)$ is a linear combination of the
$\Lambda^k$'s, $k =1, 2, 3, 4$. Lemma \ref{l:4squares}(i) implies that $A_2$ is
divisible by $\|X\|^2$. Then by \eqref{eq:plkp=rxk} with $k=2$, the matrix
$A_1 \, \mathrm{diag} \{\mu_4^2, \mu_5^2, \mu_6^2\} \, A_1^t$ is also divisible
by $\|X\|^2$, and so $\|X\|^2$ divides $A_1$, as well.

Let now $d=3$, so that all the entries of $A(X)$ are homogeneous cubic polynomials.
Let $\hat A$ be an $8 \times 5$ matrix obtained from $A(X)$ by crossing out the
last column. We can assume that $\hat A$ is not divisible by $\|X\|^2$, as otherwise
from \eqref{eq:plp=rx} it follows that $A(X)$ is also divisible by $\|X\|^2$.

As $\mu_1^{-1} < \mu_2^{-1} < \ldots < \mu_6^{-1}$, all the numbers
$\alpha_s = \mu_s^2-\mu_s\mu_6$ with $s \le 5$ have the same sign. So a
$5 \times 5$ matrix $\tilde \Lambda =\mathrm{diag} \{\alpha_1, \ldots , \alpha_5\}$
is definite. By \eqref{eq:plkp=rxk}, the polynomial matrix
$A(X) (\Lambda^2 - \mu_1 \Lambda) A(X)^t= \hat A \tilde \Lambda \hat A^t$
is divisible by $\|X\|^2$. Then by Lemma~\ref{l:4squares}(ii), we can find
polynomial vectors $U$ and $Q$ of dimension $8$ and $5$, respectively, and a
polynomial matrix $B$, such that
\begin{gather}\label{eq:rk1+x^2}
\hat A=UQ^t + \|X\|^2 B, \quad \text{and} \\
\sum\nolimits_{s=1}^5 \alpha_s q_s^2 \quad \text{is divisible by $\|X\|^2$},\label{eq:q^2}
\end{gather}
where $Q^t = (q_1, \ldots, q_5)$. Note that none of the $q_s$'s is divisible by $\|X\|^2$,
because otherwise by Lemma~\ref{l:4squares}(i), the equation \eqref{eq:q^2} implies that
all of
them are, and so $\hat A$ is divisible by $\|X\|^2$.
Moreover, from
\eqref{eq:plp=rx},
\begin{equation}\label{eq:hat a^2}
  \hat A^t \hat A = \|X\|^6 I_8,
\end{equation}
and, in particular,
\begin{equation}\label{eq:u^2}
  \sum\nolimits_{j=1}^8 u_s^2\quad \text{is divisible by $\|X\|^2$},
\end{equation}
We can choose, for every $j$ and $s$, the components $u_j, q_s$ of $U, Q$
to be of the lowest possible degree in the cosets $u_j + \bigl( \|X\|^2\bigr)$
and $q_s + \bigl( \|X\|^2\bigr)$, respectively. Then for all $s, \; q_s \ne 0$;
some of the $u_j$'s can vanish, but not all of them, as otherwise $\hat A$
is divisible by $\|X\|^2$.

Let $l_s = \deg q_s$ and $r_j= \deg u_j$ for those $j$ for which $u_j \ne 0$.
Then $l_s + r_j \le 3$, because otherwise from \eqref{eq:rk1+x^2} it will follow
that the product of the highest degree terms of $q_s$ and $u_j$ is divisible by
$\|X\|^2$, and so we can decrease the degree of one of them. Also from
\eqref{eq:rk1+x^2}, $l_s + r_j \ge 2$, as the entries of the matrix
$\hat A-\|X\|^2 B$ have no terms of degree lower than $2$. If for some
$s, j, \; l_s + r_j = 2$, then by \eqref{eq:rk1+x^2}
$q_s u_j = (\hat A)_{js}-\|X\|^2 B_{js}$, and so $q_s u_j$ is divisible by
$\|X\|^2$, as $\deg (\hat A)_{js} = 3$. Then $u_j$ is divisible by $\|X\|^2$
(as $q_s$ is not) which contradicts to the choice of $u_j$.

So for all $s$ and all $j$ with $u_j \ne 0, \; l_s + r_j = 3$, hence
$l_s = l, \; r_j = r$, and $l+r = 3$. Note that by \eqref{eq:q^2}, $r \ge 2$,
as a sum of squares of less than eight nonzero linear polynomials cannot be
divisible by $\|X\|^2$. Similarly, by \eqref{eq:u^2}, $l \ge 1$, and if $l=1$,
then all the $u_j$'s are non-zero. So, for all $s, j, \; \deg q_s = 2, \deg u_j =1$.
As $q_s u_j = (\hat A)_{js}-\|X\|^2 B_{js}$, by \eqref{eq:rk1+x^2}, none of the
$q_s, u_j$ contain lower degree terms, that is, $q_s, u_j$ are homogeneous.

So there exists an $8 \times 8$ real matrix $V$ such that $U=VX$. It follows
immediately from \eqref{eq:u^2}, that $V$ is proportional to an orthogonal matrix:
$VV^t = \beta I_8, \, \beta>0$ (if $\beta = 0$, then $U=0$, and so $\hat A$ is
divisible by $\|X\|^2$).

Consider the polynomial vector $P_1(X)$, the first column of the matrix $\hat A$.
We have
\begin{equation}\label{eq:p1}
P_1(X) = \< M_1 X, X \> VX + \|X\|^2 B_1(X),
\end{equation}
where $M_1$ is the matrix of the quadratic form $q_1$ and $B_1(X)$ is the first
column of the matrix $B(X)$. Then the equation \eqref{eq:pp} implies
$$
-X\|X\|^8 = P_1(P_1(X))= \< M_1 P_1(X), P_1(X)\> VP_1(X) + \|P_1(X)\|^2 B_1(P_1(X)).
$$
But $\|P_1(X)\|^2 = \|X\|^6$ by \eqref{eq:P_i2}, so the vector
$\< M_1 P_1(X), P_1(X)\> VP_1(X)$ must be divisible by $\|X\|^2$. Substituting $P_1(X)$
from \eqref{eq:p1} we find that $\|X\|^2$ divides $\< M_1 X, X \>^3 \<M_1VX,VX\>V^2X$.
Since the polynomial $\< M_1 X, X \>( = q_1)$ is not divisible by $\|X\|^2$, the vector
$\< M_1 VX, VX\> V^2X$ must be divisible by $\|X\|^2$. As
$V$ is a nonzero multiple of an orthogonal matrix, the
matrix $V^2$ is nonzero, and so $\< M_1 VX, VX\> = \gamma \|X\|^2$ for some
$\gamma \in \mathbb R$. Then $V^tM_1V = \gamma I_8$, so $M_1 = \gamma\beta^{-1} I_8$,
and $q_1 = \< M_1 X, X \> = \gamma\beta^{-1} \|X\|^2$ which is a contradiction.

\vskip .3cm

Thus, all the entries of the matrix $A(X)= \bigl(P_1(X)\,|\,\ldots\,|\, P_\nu(X)\bigr)$
are linear in $X$, so there exist linear operators $J_s$ in $\Ro$ such that
$P_s(X)= J_sX, \; s=1, \ldots, \nu$.
It easily follows from \eqref{eq:P_i2} that $J_s$'s are orthogonal, skew-symmetric,
and anticommute. Also, by \eqref{eq:plp=rx}, for any $X, Y \in \Ro$,
$$
R_XY = A(X) \Lambda A(X)^t Y =
\bigl(J_1X\,|\ldots|\,J_\nu X \bigr) \diag\{\mu_1, \ldots, \mu_\nu\}
\bigl(J_1X\,|\ldots|\,J_\nu X \bigr)^t Y=
\sum_{s=1}^\nu \mu_s \< J_sX, Y\> J_sX,
$$
hence the algebraic curvature tensor $R$ has a $\Cliff(\nu)$-structure.

\section{Proof of Proposition~\ref{p2}}\label{s:cliff8}

Let the curvature tensor $R$ of the manifold $M^8$ have a Clifford structure
$\Cliff(\nu)$, that is, there exist $\nu$ fields of anticommuting skew-symmetric
orthogonal operators $J_1, \dots , J_\nu$, and the functions
$\la_0, \la_1, \dots, \la_\nu \,(\la_i \ne \la_0$ for $i > 0$)
such that for any three vectors $X, Y, Z$,
\begin{multline} \label{eq:40}
R(X, Y) Z = \la_0 (\< X, Z \> Y - \< Y, Z \> X) \\
+ \sum\nolimits_{i=1}^\nu \tfrac 13 (\la_i - \la_0) (2 \< J_iX, Y \> J_iZ +
\< J_iZ, Y \> J_iX - \< J_iZ, X \> J_iY)
\end{multline}
We work in a neighbourhood $\mathcal{U}$ of a generic point on $M^8$, so
that we assume both $\nu$ and the number of pairwise nonequal $\la_i$'s to be locally
constant. As $M^8$ is Einstein, and hence analytic, it is sufficient to prove the
Proposition for $\mathcal{U}$ only.

For $\nu \le 2$ the proof follows from the results of \cite{GSV, Nhjm}, so we will
assume $3 \le \nu \le 7$. Also, by the result of \cite{Ndga}, we can assume that
the Jacobi operator has at least three different eigenvalues.

To avoid considering too many cases, we start with the following Lemma.

\begin{lemma}\label{l:form of R} One of two possibilities can occur:

\emph{(a)} There exist $7-\nu$ skew-symmetric orthogonal operators $J_{\nu+1}, \ldots , J_7$
such the operators $J_1, \ldots, J_7$ anticommute, $J_1J_2\ldots J_7 = I_8$, and the
curvature tensor has the form
\begin{equation}\label{eq:Rcase1}
  R(X, Y) Z = \sum\nolimits_{i=1}^7 \tfrac {\la_i}3
(2 \< J_iX, Y \> J_iZ + \< J_iZ, Y \> J_iX - \< J_iZ, X \> J_iY).
\end{equation}
Some of the $\la_i$'s can be equal, but up to relabelling,
we can assume that $\la_4 \ne \la_5,\, \la_6 \ne \la_7$.

\emph{(b)} $\nu = 3, \; J_1J_2 = J_3$, and the curvature tensor has the form
\begin{equation}\label{eq:Rcase2}\begin{split}
  R(X, Y) Z &= \la_0 (\< X, Z \> Y - \< Y, Z \> X) \\&+
  \sum\nolimits_{i=1}^3 \tfrac 13 (\la_i - \la_0)
(2 \< J_iX, Y \> J_iZ + \< J_iZ, Y \> J_iX - \< J_iZ, X \> J_iY).
\end{split}
\end{equation}
\end{lemma}

Note that manifolds with curvature tensor \eqref{eq:Rcase2} were studied
in \cite[Theorem 7.1]{GSV}, but our assumption
is a little weaker: all the $\la_i$'s can be nonconstant.

\begin{proof}
Let $\Cl(\nu)$ be a Clifford algebra on $\nu + 1$ generators
$1, x_1, \dots , x_\nu$. Then the map $\rho: \Cl(\nu) \to \mathbb R^n$ defined on
generators by $\rho(1) = I_n, \; \rho(x_i) = J_i$ is a representation of
$\Cl(\nu)$ in $\mathbb R^n$.

From the representation theory of Clifford algebras
\cite{ABS,H,LM},
we know that for $\nu = 4, 5, 6$, there exists a unique (up to equivalency) irreducible
representation $\rho_{\nu}$ of $\Cl(\nu)$ which is a representation in $\Ro$. The
Clifford algebra $\Cl(7)$ has two inequivalent irreducible representations, both in
$\Ro$. In fact, $J_1J_2 \ldots J_7 = \rho(x_1x_2 \ldots x_7)$ is always $\pm I_8$,
and the choice of the sign determines one of two representations $\rho_7^{\pm}$.
Up to equivalency, the representations $\rho_7^{\pm}$ can be chosen in such a way
that the operators $\tilde J_i  =\rho_7^{\pm}(x_i)$ are orthogonal and skew-symmetric.
Similarly, $\Cl(3)$ has two inequivalent irreducible representations $\rho_3^{\pm}$,
both in $\mathbb R^4$, differing by the sign of $\rho_3^{\pm}(x_1x_2x_3) = \pm I_4$.

First assume $\nu=7$. If $J_1 \ldots J_7 = -I_8$, replace $J_7$ by $-J_7$ (this does not
change the form of the curvature tensor). As for any nonzero $X \in \Ro$, the vectors
$X, J_1X, \ldots, J_7X$ form an orthogonal basis, with all the vectors of length $\|X\|$,
we have $\|X\|^2 Y - \<X, Y\> X = \< J_iX, Y \> J_iX$, and so by polarization
$\< X, Z \> Y - \< Y, Z \> X = \sum_{i=1}^7 \tfrac 13 (2 \< J_iX, Y \> J_iZ +
\< J_iZ, Y \> J_iX - \< J_iZ, X \> J_iY)$. Substituting this to \eqref{eq:40}, we obtain
\eqref{eq:Rcase1}.

Let now $\nu =4, 5, 6$. Then $\rho_{\nu}$ is equivalent to a
restriction of $\rho_7^+$ to the subalgebra $\Cl(\nu) \subset \Cl(7)$ generated by
$1, x_1, \ldots, x_\nu$. Denote $\tilde{J}_i = \rho_7^+(x_i)$ for $i = 1, \ldots , 7$.
Then there exist $T \in \mathrm{GL}(8)$ such that $J_i = T \tilde J_i T^{-1}$ for
$i =1, \ldots , \nu$. As both $J_i$ and $\tilde J_i$ are skew-symmetric, we get
$J_i (TT^t) = (TT^t)J_i$, and so every eigenspace of the symmetric operator $TT^t$
is an invariant subspace of $\rho_{\nu}$. Since $\rho_{\nu}$ is irreducible,
$TT^t$ is proportional to the identity operator, and so we can assume that
$T \in \mathrm{O}(8)$. Then for $i = \nu+1, \ldots , 7$, the operators
$J_i = T \tilde J_i T^{-1}$ are skew-symmetric and orthogonal. The anticommutativity
of $J_1, \ldots, J_7$ follows from that in $\Cl(7)$. Thus, the curvature tensor can
be written in the form \eqref{eq:40}, with the summation ranging up to $7$, and
with $\la_i=\la_0$ for $\nu < i \le 7$. We then proceed as in the case $\nu = 7$.

Finally, consider the case $\nu = 3$. We have three inequivalent
(reducible) representations of $\Cl(3)$ in $\Ro$ depending on the
symmetric operator $S = J_1J_2J_3$: either $S = \pm I_8$, or it has
two orthogonal $4$-dimensional eigenspaces $V_{\pm}$ corresponding
to eigenvalues $\pm 1$. The first case gives us the curvature
tensor of the form \eqref{eq:Rcase2} (up to switching $J_3$ to $-J_3$,
if $S = I_8$). For the representation $\rho$ in the second case,
consider the restriction of $\rho_7^+$ to the subalgebra $\Cl(3) \subset \Cl(7)$
generated by $1, x_1, x_2, x_3$. Let $\tilde J_i =\rho_7^+(x_i), \; i =1, 2, 3$.
We get a representation of $\Cl(3)$ in $\Ro$ which is equivalent to $\rho$,
as $\tilde J_1 \tilde J_2 \tilde J_3 \ne \pm I_8$. 
Then for some $T \in \mathrm{GL}(8), \; J_i = T \tilde J_i T^{-1}$, and so
$J_i (TT^t) = (TT^t)J_i, \; i = 1, 2, 3$. Hence the symmetric positively defined
operator $TT^t$ has two eigenvalues, $\alpha_{\pm}^2$, with eigenspaces $V_{\pm}$,
respectively. Define an operator $\tilde T$ by
$\tilde T_{|V_{\pm}} = \alpha_{\pm}^{-1} T_{|V_{\pm}}$. Then $\tilde T$ is orthogonal,
and for $i =1, 2, 3, \; J_i = \tilde T \tilde J_i {\tilde T}^{-1}$, as $V_{\pm}$
are the invariant subspaces of $J_1, J_2, J_3$. We then define $J_i,\; 4 \le i \le 7$
by $J_i = \tilde T \,\rho_7^+(x_i) \,{\tilde T}^{-1}$ and proceed as in the case
$\nu = 4,5,6$.

The eigenvalues of the Jacobi operator of the curvature tensor \eqref{eq:Rcase1}
are $\la_1, \ldots, \la_7$. As at least three of them are pairwise nonequal, they
can be relabelled to satisfy $\la_4 \ne \la_5,\; \la_6 \ne \la_7$ (this can be
done by an even permutation, not to violate the condition $J_1\ldots J_7 =I_8$).
\end{proof}

We are going to show that in each of the cases (a) and (b) of Lemma~\ref{l:form of R},
the manifold $M^8$ is rank-one symmetric.

\subsection{Case (a)} \label{ss:casea}

In this case, the manifold carries a Hermitian almost octonion structure
defined by the operators $J_1, \ldots, J_7$. The $\Cl(7)$-module
$\Ro$ for the representation $\rho_7^+$ can be described in terms of the
octonion algebra $\mathbb O$ as follows. Let $1, e_1, \ldots, e_7$ be
a fixed \textit{canonical basis} of generators for $\mathbb O$, so that
$e_i^2 = -1, \, e_ie_j+e_j e_i = 0, \; i \ne j$ (we will use the canonical
basis with multiplication table given in \cite[Ch.~3G]{Besse}). Then there
exists an orthogonal map $\psi: \Ro \to \mathbb O$ such that for every
$i = 1, \dots, 7$ and every $X \in \Ro, \; \psi(J_iX) = \psi(X) e_i$.

\begin{remark}
Explicitly, one can construct $\psi$ from the given $J_i$'s as follows:
consider the set $S$ of triples $(i, j, k)$ with $e_ie_j = e_k$, according
to the multiplication table. Any two triples from $S$ either contain
exactly one element in common, or can be obtained from one another by a
cyclic permutation. Then the operators $\{J_iJ_jJ_k,\, (i,j,k) \in S\}$ in
$\Ro$ are orthogonal, symmetric and pairwise commute, hence they have the
same set of eigenvectors. One can show that it is possible to chose a
common eigenvector $X^0$ whose eigenvalue is $+1$ for all the above
operators. Then $\psi$ is defined by $\psi(X^0) =1, \; \psi(J_iX^0)=e_i$.
\end{remark}

With some abuse of language, we identify, by $\psi$, vectors from $\Ro$
with octonions, so that a vector $X \in \Ro$ with coordinates
$(X_0, X_1, \ldots , X_7)$ is the same as the octonion
$X_01+X_1e_1+\ldots +X_7e_7$. It is convenient to introduce the following
notation:
$$
J_uX := \sum\nolimits_{i=1}^7 u_i J_iX = X u,
$$
for $X \in \Ro = \mathbb O, \; u = \sum_{i=1}^7 u_i e_i \in \mathbb O' := 1^\perp$.
We will also abbreviate $J_{e_i}$ to $J_i$.

Throughout this Section we will use, without explicitly referring, the following
identities (see \cite[Sec.~2]{BG} for details):
\begin{gather*}
a^*:=2 \<a,1\> 1 - a,\quad (ab)^* = b^*a^*,\quad a^{**} = a, \quad aa^* = \|a\|^21,
\quad \|ab\| = \|a\|\|b\|,\\
\<a^*, b^*\> = \<a,b\>, \quad \<a, bc\> = \<b^*a, c\> = \<ac^*, b\>,\\
a(ab) = (a^2) b, \quad (ba)a = b (a^2), \quad a(a^*b) = (ba)a^* = \|a\|^2 b,\\
(ab^*)c+(ac^*)b = 2\<b, c\>a, \quad a(b^*c)+ b(a^*c) = 2\<a, b\>c,
\end{gather*}
for any $a, b, c \in \mathbb{O}$. Note that $1^* = 1, \, e_i^* = -e_i$.
The algebra $\mathbb{O}$ is a \textit{division algebra}, in particular, any nonzero
octonion is invertible, with $a^{-1} = \|a\|^{-2} a^*$.

We will also use \textit{bioctonions} $\mathbb O \otimes \mathbb C$, the algebra over
the $\mathbb{C}$ with the same canonical basis of generator as $\mathbb O$. As all
the above identities are polynomial, they still hold for bioctonions, with Euclidean
scalar product in $\mathbb{C}^8$, the underlying linear space of
$\mathbb O \otimes \mathbb C$. However, bioctonion algebra is no longer a division
algebra (for example, because it has zero-divisors: $(i1+e_1)(i1-e_1)=0$).

The following Lemma will be in use later in the proof.
\begin{lemma} \label{l:octonions}
Let $e \perp 1$ be a unit octonion, and let $L, F: \mathbb{O} \to \mathbb{O}$ be
\emph{(}$\mathbb R$-\emph{)}linear operators satisfying
\begin{equation} \label{eq:YeL=YF}
(Ye)L(Y)=YF(Y)\quad \text {for all} \; Y \in \mathbb O.
\end{equation}
Then there exist $a, b, p \in \mathbb O$ such that
$$
L(Y)= \<a,Y\>\,1 + \<b, Y\>\,e + Y^*p.
$$
\end{lemma}

\begin{proof}
The set of linear operators $L$ for which there exists $F$ such that $L, F$ satisfy
\eqref{eq:YeL=YF} is a linear subspace $\mathcal S \subset \mathrm{End}(\Ro)$. For any $p \in \mathbb O$,
an operator $L(Y) = Y^*p$ is in $\mathcal S$, as
$(Ye)(Y^*p) = 2\<Ye,Y\>p - Y((Ye)^*p) = Y(-(Ye)^*p)$, and we can take $F(Y) = -(Ye)^*p$.

So, given an operator $L$ satisfying \eqref{eq:YeL=YF} with some $F$, we can replace it by
an operator $\hat L$ defined as $\hat L(Y) = L(Y) - Y^*L(1)$, still satisfying \eqref{eq:YeL=YF}
and having the property $\hat L(1) = 0$. The operator $\hat F$ corresponding to $\hat L$
by \eqref{eq:YeL=YF} must also vanish at $Y=1$.

Denote $\mathcal L =\left(\Span (1,e)\right)^\perp$. Let $u \in \mathcal L$ be a nonzero octonion
and $t$ a real number. Substituting $Y = t1+u$ in \eqref{eq:YeL=YF} we get:
$(te +ue) \hat L(u) = (t1+u) \hat F(u)$, so $\hat F(u) = e \hat L(u)$ and
$(ue)\hat L(u)=u\hat F(u)$, hence
$$
(ue)\hat L(u)=u(e \hat L(u)).
$$
It follows that for all $u \in \mathcal L,\; u \ne 0$, the octonion $\hat L(u)$ lies in an associative
subalgebra of $\mathbb O$ generated by $u, e$. This subalgebra is $\Span (1,e,u,ue)$ and is isomorphic to
$\mathbb H$. Then for some functions
$\alpha,\beta,\gamma,\delta: \mathcal L \to \mathbb R$ we have
\begin{equation}\label{eq:hatL}
\hat L(u) = \alpha(u)1+\beta(u)\,e+\gamma(u)\,u+\delta(u)\,ue.
\end{equation}
The functions $\alpha(u) = \<\hat L(u),1\>$ and $\beta(u) = \<\hat L(u),e\>$ are linear,
so there exist octonions $\hat a, \hat b$ such that $\alpha(u) = \<\hat a,u\>,\; \beta(u)=\<\hat b,u\>$.
We want to show that $\gamma$ and $\delta$ are constants. Indeed, from \eqref{eq:hatL} we get
$\tilde L(u) = \gamma(u)u+\delta(u)ue$, where $\tilde L$ is a linear operator defined by
$\tilde L(u) = \hat L(u) -\<\hat a,u\>1-\<\hat b,u\>e$. Then
$\gamma(u) = \|u\|^{-2} C(u),\; \delta(u) = \|u\|^{-2} D(u)$, where
$C(u) = \<\tilde L(u), u\>,\; D(u) = \<\tilde L(u), ue\>$ are quadratic forms on $\mathcal L$, and
$\|u\|^2 \tilde L(u) = C(u)\, u + D(u)\, ue$. Taking the squared norm of the both sides we obtain
$\|u\|^2 \|\tilde L(u)\|^2 = C^2(u) + D^2(u)$. As the sum of squares of polynomials $C, D$ in six
real variables is divisible by $\|u\|^2$, each of them must be divisible by $\|u\|^2$ (see the
Remark in Subsection~\ref{ss:ca}). So
$\gamma = \|u\|^{-2} C(u)$ and $\delta = \|u\|^{-2} D(u)$ are constants. Back to \eqref{eq:hatL}, we get
$$
\hat L(u) = \<\hat a, u\>\,1+\<\hat b,u\>\,e+\gamma \,u+ \delta \,ue.
$$
Then for an arbitrary $Y \in \mathbb O$,
\begin{align*}
\hat L(Y) &= \hat L(\<Y,1\>1 + \<Y,e\>e + (Y-\<Y,1\>1 - \<Y,e\>e))\\
&=\<Y,e\> \hat L(e) + \<\hat a, Y-\<Y,1\>1 - \<Y,e\>e\>1+\<\hat b,Y-\<Y,1\>1 - \<Y,e\>e\>e\\
&+ \gamma (Y-\<Y,1\>1 - \<Y,e\>e)+ \delta (Y-\<Y,1\>1 - \<Y,e\>e)e \\
&=\<Y,e\> \hat L(e) + \<a,Y\>1 + \<b, Y\> e + Y^*p,
\end{align*}
where $a = \hat a - \<\hat a, 1\>1 - \<\hat a, e\>e +\gamma 1 + \delta e,\;
b = \hat b - \<\hat b, 1\>1 - \<\hat b, e\>e  - \gamma e + \delta 1,\;
p = -\gamma 1 - \delta e$. To finish the proof, it remains to show that
the octonion $\hat L(e)$ lies in $\Span(1, e)$. To see that, substitute the above
expression for $\hat L(Y)$ to \eqref{eq:YeL=YF}. Then
$$
\<Y,e\> (Ye)\hat L(e)= Y \tilde F(Y),
$$
with a linear operator $\tilde F(Y)= \hat F(Y) - (\<a, Y\> e - \<b, Y\> 1 - (Ye)^*p)$.
Multiplying the above equation by $Y^*$ from the left we get
$\<Y,e\> Y^*((Ye)\hat L(e))= \|Y\|^2 \tilde F(Y)$,
so the octonion $Y^*((Ye)\hat L(e))$ viewed as an $8$-dimensional vector, polynomial in $Y$, is
divisible by $\|Y\|^2$. Then there exists an octonion $f$ such that $Y^*((Ye)\hat L(e)) = \|Y\|^2f$
which implies $(Ye)\hat L(e) = Yf$, for all $Y$. If $Y=1$, this gives $f = e \hat L(e)$, and so
$(Ye)\hat L(e) = Y(e \hat L(e))$. Therefore $\hat L(e)$ belongs to all associative subalgebras of
$\mathbb O$ generated by $e, Y$, with $Y$ an arbitrary octonion. This is only possible when
$\hat L(e) \in \Span(1, e)$.
\end{proof}

In the next Lemma, we calculate covariant derivatives of the operators $J_i$, for later
use in the second Bianchi identity.

\begin{lemma}\label{l:nablaJ}
There exist linear functionals $A_i, B_{ij}$ on $\Ro, \; i,j =1, \ldots, 7$, with
$B_{ij}=-B_{ji}$ such that for all $U \in \Ro$
\begin{equation}\label{eq:nablaJ}
  \n_UJ_i = \sum\nolimits_{j\ne i} B_{ij}(U)J_j + \sum\nolimits_{j\ne i} A_j(U)J_iJ_j
  = \sum\nolimits_j B_{ij}(U) J_j + J_i J_{A(U)} + A_i (U) I_8,
\end{equation}
where $A: \mathbb O \to \mathbb O' \, (= \mathbb O \cap 1^\perp)$ is a linear operator
defined by $A(U) = \sum_j A_j(U) e_j$.

\end{lemma}
\begin{proof} The operator $\n_UJ_i$ depends linearly on $U$ and is skew-symmetric.
As the kernel of $\rho_7^+$ is the ideal of $\Cl(7)$ generated by
$x_1 \ldots x_7 -1$, the $28$ skew-symmetric operators
$\{J_j, J_jJ_k\},\; 1 \le j < k \le 7$
are linearly independent. In fact, they form a basis for $\mathfrak o(8)$, the space
of skew-symmetric operators in $\Ro$, as $\dim \mathfrak o(8) = 28$.
So $\n_UJ_i = \sum_j B_{ij}(U)J_j + \sum_{k,j} C_{i;kj}(U)J_kJ_j$ for some linear
functionals $C_{i;kj}, B_{ij}$, with $C_{i;kj}= - C_{i;jk}$.

As $\n_U(J_i^2) = 0$, the operator $J_i\n_UJ_i$ is skew-symmetric, so its symmetric part
$-B_{ii}(U) I_8 + \sum_{k,j \ne i} C_{i;kj}(U)J_iJ_kJ_j$
vanishes. Acting by $J_i$ from the left, we get a linear combination of $J_i, J_kJ_j$
which is zero, therefore for all $i$ and for all $j, k \ne i, \quad B_{ii} = C_{i;kj} = 0$
which gives
$\n_UJ_i = \sum_{j\ne i} B_{ij}(U)J_j + \sum_{j \ne i} C_{i;j}(U)J_iJ_j$, where
$C_{i;j} = 2 C_{i;ij}$.

The equation $\n_U(J_iJ_s + J_sJ_i) = 0$ implies that the operator
$J_i \n_UJ_s + J_s \n_UJ_i = \sum_{j\ne s} B_{sj}(U)J_iJ_j + \sum_{j \ne s}
C_{s;j}(U)J_iJ_sJ_j
+ \sum_{j\ne i} B_{ij}(U)J_sJ_j + \sum_{j \ne i} C_{i;j}(U)J_sJ_iJ_j$ is skew-symmetric,
for all $s \ne i$. Extracting the symmetric part we get:
$-(B_{si}(U) + B_{is}(U)) I_8 + \sum_{j \ne i,s} (C_{s;j}(U)- C_{i;j}(U))J_iJ_sJ_j = 0$.
So
$B_{si}(U) + B_{is}(U) = 0$ and $C_{i;j}(U)$ depends only on $j$ for $j \ne i$. Denoting
$A_i (U) =C_{i;j}(U), \; j \ne i$ we obtain \eqref{eq:nablaJ}.
\end{proof}

We will use the second Bianchi identity
\begin{equation}\label{eq:Bi}
(\n_UR)(X, Y, Y, X) + (\n_XR)(Y, U, Y, X) + (\n_YR)(U, X, Y, X) = 0.
\end{equation}

From \eqref{eq:Rcase1}, $R(V, Z) V = \sum_{i=1}^7 \la_i \< J_iV, Z \> J_iV$,
and so for any three vectors $V, Z, W$,
\begin{equation}\label{eq:nablaR}
\begin{split}
  (\n_W R)(V, Z) V &= \sum\nolimits_{i=1}^7 W(\la_i) \< J_iV, Z \> J_iV \\&+
\sum\nolimits_{i=1}^7 \la_i (\< (\n_WJ_i)V, Z \> J_iV + \< J_iV, Z \> (\n_WJ_i)V).
\end{split}
\end{equation}

Substituting this to \eqref{eq:Bi} we obtain
\begin{equation}\label{eq:old7}
\begin{split}
&\sum\nolimits_{i=1}^7 (X(\la_i) \< J_iY, U \> \< J_iY, X \> +
Y(\la_i) \< J_iU, X \> \< J_iY, X \>
- U(\la_i) \< J_iY, X \>^2 ) \\
+ &\sum\nolimits_{i=1}^7 \la_i \, \< J_iY, X \>
(2 \< (\n_UJ_i)X, Y \> + \< (\n_XJ_i)Y, U \> + \< (\n_YJ_i)U, X \>)
\\-&\sum\nolimits_{i=1}^7 \la_i \, (\< J_iY, U \> \< (\n_XJ_i)X, Y \>
+ \< J_iX, U \> \< (\n_YJ_i)Y, X \>) = 0.
\end{split}
\end{equation}

\begin{lemma}\label{l:a(U)}
The linear operator $A$ introduced in Lemma~\ref{l:nablaJ} has the form
\begin{equation}\label{eq:A(U)}
A(U) = U^*m - \<U, m\> \, 1, \quad \text{for some} \; m \in \mathbb O.
\end{equation}
\end{lemma}

\begin{proof}
The equation \eqref{eq:old7} is a polynomial equation in $24$ real variables, the
coordinates of the vectors $X, Y, U$. It must still hold if we allow $X, Y, U$ to
be complex and extend the operators $J_i$ to $\mathbb C^8$ by complex linearity.
Let $Y \in \mathbb C^8$ be a nonzero isotropic vector (bioctonion): $\|Y\|^2 = 0$.
Denote $\mathcal JY = Y (\mathbb O \otimes \mathbb C) =
\Span_{\mathbb{C}}(J_1Y, \ldots, J_7Y)$. The space $\mathcal JY$ is isotropic:
the scalar product of any two vectors from $\mathcal JY$ is zero. Choosing
bioctonions $X, U \in \mathcal JY$, we get $\< J_iY, X \> = \< J_iY, U \> = 0$,
for all $i =1, \ldots 7$. Substituting this to \eqref{eq:old7} we obtain
$\sum_{i=1}^7 \la_i \<(\n_YJ_i)Y, X \> \<J_i X, U\> =0$ which simplifies to
$\sum_{i=1}^7 \la_i \<J_iJ_{A(Y)}Y, X \> \<J_i X, U\> =0$ by
\eqref{eq:nablaJ}. So $J_{A(Y)}Y \perp \sum_{i=1}^7 \la_i  \<J_i X, U\> J_iX$
for all $X, U \in \mathcal JY$.
From this and the fact that $J_{A(Y)}Y \in \mathcal JY \subset (\mathcal JY)^\perp$
we get:
\begin{equation}\label{eq:JaYperp}
J_{A(Y)}Y \perp \mathcal LY := \mathcal JY +
\Span_{\mathbb{C}}\left\{\sum\nolimits_{i=1}^7 \la_i  \<J_i X, U\> J_iX, \;|\; X, U \in \mathcal JY\right\}
\end{equation}

The linear space $\mathcal LY$ is at
most $7$-dimensional, as $Y \perp \mathcal LY$. If we can produce at least one
isotropic vector $Y$ such that $\dim \mathcal LY = 7$, then it will follow
from \eqref{eq:JaYperp} that for an open set of $Y$'s on the isotropic cone
$\|Y\|^2 = 0$, the vector $J_{A(Y)}Y$ is parallel to $Y$. As $J_{A(Y)}Y$ depends
algebraically (in fact, quadratically) of $Y$, we will have $J_{A(Y)}Y \parallel Y$
for all isotropic vectors $Y$.

Take $Y=i1+e_1$. Using the multiplication table of \cite[Ch.3G]{Besse} we find
$$
\mathcal J(i1+e_1) = \Span_{\mathbb{C}}(i1+e_1, ie_2 + e_3, ie_4+e_5, ie_6-e_7).
$$
Evaluating the expression $\sum_{i=1}^7 \la_i  \<J_i X, U\> J_iX$ with
$X = ie_2 + e_3, U = ie_6-e_7$, with $X = ie_2 + e_3, U = ie_4+e_5$, and with
$X = ie_4+e_5, U = ie_2 + e_3$, we get, respectively,
$-2(\lambda_4-\lambda_5) (ie_6 +e_7), \, 2(\lambda_6-\lambda_7) (-ie_4 +e_5), \,
-2(\lambda_6-\lambda_7) (ie_2 -e_3)$. As
$\lambda_4 \ne \lambda_5, \, \lambda_6 \ne \lambda_7$ by the assumption
(Lemma~\ref{l:form of R}), the vectors
$ie_6 +e_7, \, -ie_4 +e_5, \, ie_2 -e_3$ lie in $\mathcal L(i1+e_1)$. It follows that
$$
\mathcal L(i1+e_1) \supset \Span_{\mathbb{C}}\{ i1+e_1, e_2, e_3, e_4, e_5, e_6, e_7 \},
$$
and so $\dim \mathcal L(i1+e_1) = 7$.

From the above, $J_{A(Y)}Y \parallel Y$ for all $Y$ with $\|Y\|^2 = 0$.
So for all $i, j$, the coordinates of the vector $J_{A(Y)}Y$ viewed
as complex polynomials of $Y_0, Y_1, \ldots , Y_7$ satisfy
$(J_{A(Y)}Y)_i Y_j = (J_{A(Y)}Y)_j Y_j \mod \left(\|Y\|^2\right)$. As the ring
$\mathbf{R}\otimes\mathbb{C} =\mathbb C[Y_0, \ldots, Y_7]/ ( \|Y\|^2 )$ is a UFD
(Fact~\ref{fact:1}, Subsection~\ref{ss:ca}),
there exists an element $r \in \mathbf{R}\otimes\mathbb{C}$ such that $\pi((J_{A(Y)}Y)_j) = r \pi(Y_j)$
for all $j$, with $\pi:\mathbb C[Y_0, \ldots, Y_7] \to \mathbf{R}\otimes\mathbb{C}$ the natural
projection. Lifting up, we find a polynomial function $f \in \pi^{-1}(r)$ and a polynomial
vector $F$ such that $J_{A(Y)}Y = f(Y) Y + \|Y\|^2 F(Y)$ for all $Y \in \mathbb C^8$.
As $A(Y)$ is a linear operator, the left-hand side is a polynomial vector with
all the components of degree $2$. As the polynomial $\|Y\|^2$ is
irreducible, we can choose $f \in \pi^{-1}(r)$ to be linear, and $F$ to be a constant.
So there exist $m, c \in \mathbb O \otimes \mathbb C$
such that for all $Y \in \mathbb O \otimes \mathbb C$,
$$
Y A(Y)\,( = J_{A(Y)}Y) = \<c, Y \> \, Y + \|Y\|^2 \, m.
$$
Multiplying by $Y^*$ from the left and dividing through by $\|Y\|^2$ we get
$A(Y) = \<c, Y\> \,1 + Y^*\,m$. Since $A(Y) \perp 1$, we must have
$\<c, Y\>+ \<m, Y\> = 0$, and so $c=-m$. Also,
as $A(Y)$ is real (lies in $\mathbb O \subset \mathbb O \otimes \mathbb C$), when $Y$
is real, $m$ has to be real, as well. So, for some $m \in \mathbb O$,
$$
A(Y) = - \<m, Y \> \, 1 + Y^* m
$$
\end{proof}

\begin{lemma}\label{l:Bandlambda}
All the $\la_k$'s are constant:
\begin{equation}\label{eq:lak=const}
\n \la_k = 0, \quad\text{for all}\quad k = 1, \ldots, 7,
\end{equation}
and the operators $B_{ij}$ introduced in Lemma~\ref{l:nablaJ} satisfy
\begin{equation}\label{eq:B_ij}
(\la_i - \la_j) B_{ij}(Y) = - (\la_i - \la_j) \<(me_j)e_i,Y\>, \quad j \ne i
\end{equation}
\end{lemma}
\begin{proof}
Substitute $X = J_kY = Ye_k$ in \eqref{eq:old7}. We have
$\<J_iY,X\> = \|Y\|^2 \delta_{ik}$, and the second sum is non-zero only when $i=k$. Moreover,
$\< (\n_UJ_k)X, Y \> = - \< J_kY, (\n_UJ_k)Y \> = 0$. Also, from
\eqref{eq:A(U)}, \eqref{eq:nablaJ} we get:
\begin{align*}
\< (\n_XJ_k)Y, U \> + \< (\n_YJ_k)U, X \> &=
\sum\nolimits_j (B_{kj}(X)\<J_jY,U\> - B_{kj}(Y)\<J_jX,U\>) \\&+
\<J_k J_{A(X)} Y + \<A(X),e_k\> Y - J_k J_{A(Y)} X - \<A(Y),e_k\> X, U\>\\&
= \sum\nolimits_j (B_{kj}(X)\<Ye_j,U\> - B_{kj}(Y)\<Xe_j,U\>) \\&+
\<(Y(X^*m) - X(Y^*m))e_k , U\> -2 (\<X,m\>\<X, U\>  + \<Y,m\>\< Y , U\>).
\end{align*}
As $(Y(X^*m) - X(Y^*m))e_k =2 (Y(X^*m))e_k =-2 (Y(X^*m))e_k^*=
-4 \<X^*m,e_k\>Y + 2 Ye_k (X^*m)^* =
-4 \<m,Xe_k\>Y + 2 X (-X^*m + 2\<X^*m,1\>1) =
4(\<Y,m\>Y + \<X,m\> X) -2 \|Y\|^2 m$, we find
\begin{multline*}
\sum\nolimits_i \la_i \, \< J_iY, X \>
(2 \< (\n_UJ_i)X, Y \> + \< (\n_XJ_i)Y, U \> + \< (\n_YJ_i)U, X \>)  \\=
\la_k \|Y\|^2 \bigl(
\sum\nolimits_j (B_{kj}(X)\<Ye_j,U\> - B_{kj}(Y)\<Xe_j,U\>)
-2\|Y\|^2 \<m,U\> + 2 (\<X,m\>\<X, U\>  + \<Y,m\>\< Y , U\>)
\bigr).
\end{multline*}
In a similar way, the third sum in \eqref{eq:old7} simplifies to
\begin{multline*}
-\|Y\|^2
\sum\nolimits_i \la_i \,( - B_{ik}(X)\<Ye_i,U\> + B_{ik}(Y)\<Xe_i,U\>
+ \<me_i,Y\>\<Ye_i, U\>  + \<me_i,X\>\<Xe_i, U\> ) \\
- \la_k \|Y\|^2 (\<X, m\>\<X, U\>  + \<Y,m\>\<Y, U\>).
\end{multline*}
Substituting all this to \eqref{eq:old7}, dividing through by $\|Y\|^2$ and
observing that the equation is linear with respect to $U$, we come to the
following octonion equation:
\begin{multline*}
X(\la_k) X + Y(\la_k) Y - \|Y\|^2 \n \la_k - 2 \la_k \|Y\|^2 m
+ \la_k (\<X, m\> X   + \<Y,m\> Y) \\-
\sum\nolimits_i \la_i (\<me_i,Y\> Ye_i  + \<me_i,X\> Xe_i)
+ \sum\nolimits_j (\la_k - \la_j)(B_{kj}(X) Ye_j - B_{kj}(Y) Xe_j) = 0,
\end{multline*}
where $X = Ye_k$. As $2\|Y\|^2 m = \|Y\|^2 m +\|X\|^2 m =
\sum_i \<m, Ye_i\>Ye_i + \<m, Y\> Y + \sum_i \<m, Xe_i\>Xe_i + \<m, X\> X
=-\sum_i (\<me_i, Y\>Ye_i + \<me_i, X\> Xe_i) +\<m, Y\> Y + \<m, X\> X$,
we get
\begin{equation*}
X(\la_k) X + Y(\la_k) Y - \|Y\|^2 \n \la_k
+ \sum\nolimits_j (\la_k - \la_j) \bigl(
(B_{kj}(X) + \<me_j,Y\>) Ye_j + (- B_{kj}(Y)+ \<me_j,X\>) Xe_j \bigr) = 0.
\end{equation*}
As $\|Y\|^2 \n \la_k = Y(Y^* \n \la_k)$ and $X = Ye_k$, the above equation can be written in the form
$$
(Ye_k)L_k(Y) = Y F_k(Y),
$$
where, for every $k =1, \ldots, 7, \; L_k$ and $F_k$ are linear operators in $\mathbb O$,
with $L_k(Y) \perp 1$, given by
\begin{equation}\label{eq:LkFk}
\begin{split}
L_k(Y)&= \sum\nolimits_j  (\la_k - \la_j) (B_{kj}(Y) + \<(me_j)e_k,Y\>) e_j,
\\
F_k(Y)&= (Ye_k)(\la_k) e_k + Y(\la_k) 1 - Y^* \n \la_k
+ \sum\nolimits_j (\la_k - \la_j) (B_{kj}(Ye_k) + \<me_j,Y\>) e_j
\end{split}
\end{equation}
Applying Lemma~\ref{l:octonions}, we find
$L_k(Y)= \<a_k, Y\>1 + \<b_k,Y\>e_k + Y^*p_k$, for some $a_k, b_k, p_k \in \mathbb O$.
Since $L_k(Y) \perp 1$, we get $a_k = - p_k$, and so
\begin{equation}\label{eq:Lk}
L_k(Y)= - \<p_k, Y\>1 + \<b_k,Y\>e_k + Y^*p_k, \qquad k =1, \ldots, 7.
\end{equation}
From \eqref{eq:LkFk}, \eqref{eq:Lk} we obtain, for $j \ne k$:
$$
\<L_k(Y), e_j \> = (\la_k - \la_j) (B_{kj}(Y) + \<(me_j)e_k,Y\>) = - \<p_ke_j, Y\>
$$
But the middle expression is symmetric with respect to $j, k$, as $B_{kj}=-B_{jk}$,
so $p_ke_j = p_je_k$, for all $j \ne k$. It follows that $p_ke_k = - p_je_j$,
hence all the $p_k$'s vanish. Again, from (\ref{eq:LkFk}, \ref{eq:Lk}) it
follows that $\<L_k(Y), e_k \> = 0 = \<b_k,Y\>$, so $b_k = 0$, and therefore
$L_k(Y)=0$. Then also $F_k(Y)=0$, and by \eqref{eq:LkFk} we obtain
\begin{gather*}
(\la_k - \la_j) B_{kj}(Y) = - (\la_k - \la_j) \<(me_j)e_k,Y\>, \quad j \ne k
\\
(Ye_k)(\la_k) e_k + Y(\la_k) 1 - Y^* \n \la_k = 0.
\end{gather*}
The former equation proves \eqref{eq:B_ij}, the latter one, after multiplying by
$Y$ from the left, implies that
$\n\la_k$ lies in the two-dimensional space $\Span(Y, Ye_k)$, for any nonzero
$Y \in \mathbb O$, and so $\n \la_k = 0$ for all $k = 1, \ldots, 7$.
\end{proof}
It now follows from \eqref{eq:nablaR} and \eqref{eq:lak=const} that for any $X, Y$,
\begin{equation*}
  (\n_X R)(X, Y) X = \sum\nolimits_i \la_i (\< (\n_XJ_i)X, Y \> J_iX + \< J_iX, Y \> (\n_XJ_i)X).
\end{equation*}
As by (\ref{eq:nablaJ}, \ref{eq:A(U)}), \;
$(\n_XJ_i)X = \sum_jB_{ij}(X)J_jX +\|X\|^2J_im - \<X,m\>J_iX+\<m,J_iX\>X$, we get
\begin{equation}\label{eq:nablaRwith}
\begin{split}
\negmedspace
(\n_X R)(X, Y) X &=
\sum\nolimits_i \la_i \left( \< \sum\nolimits_jB_{ij}(X)J_jX +\|X\|^2J_im - \<X,m\>J_iX+\<m,J_iX\>X, Y \> J_iX \right.\\
&\left.+ \< J_iX, Y \> (\sum\nolimits_jB_{ij}(X)J_jX +\|X\|^2J_im - \<X,m\>J_iX+\<m,J_iX\>X) \right)\\
&=
\sum\nolimits_{ij} (\la_i-\la_j) B_{ij}(X)\<J_jX, Y \> J_iX- 2\<X,m\>\sum\nolimits_i \la_i \<J_iX,Y\>J_iX  \\
&+\sum\nolimits_i \la_i \left(\< \|X\|^2J_im +\<m,J_iX\>X, Y \> J_iX + \< J_iX, Y \> (\|X\|^2J_im +\<m,J_iX\>X)\right).
\end{split}
\end{equation}
Now from \eqref{eq:B_ij}
\begin{align*}
&\sum\nolimits_{ij} (\la_i-\la_j) B_{ij}(X)\<J_jX, Y \> J_iX=
-\sum\nolimits_{ij} (\la_i-\la_j) \<(me_j)e_i,X\> \<J_jX, Y \> J_iX\\
&=\sum\nolimits_{ij} \la_j \<(me_j)e_i,X\> \<J_jX, Y \> J_iX
-\sum\nolimits_{ij} \la_i \<(me_j)e_i,X\> \<J_jX, Y \> J_iX\\
&=-\sum\nolimits_i \left\langle\sum\nolimits_j \la_j \<J_jX,Y\>J_jm, J_iX\right\rangle J_iX\\
&\hphantom{=\;}-\sum\nolimits_i \la_i\left(\sum\nolimits_{j\ne i}\<J_im, J_jX\>\<J_jX,Y\> - \<m,X\>\<J_iX,Y\>
\right)J_iX.
\end{align*}
As $\sum_i \<U,J_iX\>J_iX = \|X\|^2 U - \<U,X\> X,\;
\sum_i \<U,J_iX\>\<V, J_iX\> = \|X\|^2 \<U, V\> - \<U,X\> \<V,X\>$, for any $U, V$, the above
expression simplifies to
\begin{multline*}
\sum\nolimits_{ij} (\la_i-\la_j) B_{ij}(X)\<J_jX, Y \> J_iX
=-\|X\|^2 \sum\nolimits_j \la_j \<J_jX,Y\>J_jm + \sum\nolimits_j \la_j \<J_jX,Y\> \<J_jm,X\>X\\
-\sum\nolimits_i \la_i (\|X\|^2 \<J_im, Y\> - 2\<m,X\> \<J_iX,Y\> - \<J_im, X\> \<Y, X\>)J_iX.
\end{multline*}
Substituting to \eqref{eq:nablaRwith} we obtain
$$
(\n_X R)(X, Y) X = 0.
$$
So $\n R = 0$ and the manifold $M^8$ is locally symmetric.
The fact that $M^8$ is rank-one symmetric follows from
\cite[Lemma 2.3]{GSV}.

\subsection{Case (b)} \label{ss:caseb}

In this case, the manifold carries an almost hypercomplex Hermitian structure
defined by the operators $J_1, J_2, J_3$, with $J_1J_2 = J_3$, and
with the curvature tensor given by \eqref{eq:Rcase2}.

The second Bianchi identity
$(\n_UR)(X, Y, Y, X) + (\n_XR)(Y, U, Y, X) + (\n_YR)(U, X, Y, X) = 0$ takes the form
\begin{equation}\label{eq:old7for3}
\begin{split}
&U(\la_0) (\< X, Y \>^2 - \|X\|^2 \|Y\|^2) \\
+ &X(\la_0) (\|Y\|^2 \< U, X \> - \< U, Y \> \< X, Y \>) + Y(\la_0) (\|X\|^2 \< U, Y \> -
\< U, X \> \< X, Y \>)   \\
+ &\sum\nolimits_{i=1}^3 (X(\la_i-\la_0) \< J_iY, U \> \< J_iY, X \> +
Y(\la_i-\la_0) \< J_iU, X \> \< J_iY, X \>
- U(\la_i-\la_0) \< J_iY, X \>^2 ) \\
+ &\sum\nolimits_{i=1}^3 (\la_i-\la_0) \, \bigl(
(2 \< (\n_UJ_i)X, Y \> + \< (\n_XJ_i)Y, U \> +
\< (\n_YJ_i)U, X \>)  \< J_iY, X \>
\\
& \hphantom{\sum\nolimits_{i=1}^3 (\la_i-\la_0) (}
- \< J_iY, U \> \< (\n_XJ_i)X, Y \>
- \< J_iX, U \> \< (\n_YJ_i)Y, X \> \bigr) = 0.
\end{split}
\end{equation}

For a nonzero vector $Y$, let $\mathcal LY = \Span(Y, J_1Y, J_2Y, J_3Y)$. Note that
$\mathcal L(J_iY) =\mathcal LY$, and that if $X \perp \mathcal LY$, then also
$\mathcal LX \perp \mathcal LY$.

Taking $U, X \perp \mathcal LY$, with $J_1X = U$, in \eqref{eq:old7for3} we obtain
\begin{equation}\label{eq:nablal0}
\|Y\|^2 U (\la_0) - (\la_1-\la_0) \< (\n_YJ_1)Y, J_1U \>  = 0
\quad\text{for all $U, Y$ with $U \perp \mathcal LY$}.
\end{equation}
On the other hand, taking $X = J_1Y,\; U \perp \mathcal LY= \mathcal LX$
in \eqref{eq:old7for3} gives
\begin{equation}\label{eq:nablal1}
 \|Y\|^2 U(\la_1) -(\la_1-\la_0) (2 \< (\n_UJ_1)X, Y \> + \< (\n_XJ_1)Y, U \> + \< (\n_YJ_1)U, X \>) = 0.
\end{equation}

Using \eqref{eq:nablal0} and the fact that the operator $J_1$ is orthogonal and skew-symmetric we get:
\begin{align*}
&\< (\n_UJ_1)X, Y \> = - \< (\n_UJ_1)Y, J_1Y \> = 0, \\
&(\la_1-\la_0)\< (\n_YJ_1)U, X \> = -(\la_1-\la_0)\< U, (\n_YJ_1)J_1Y \> = -(\la_1-\la_0)\< J_1U, (\n_YJ_1)Y \> =
-\|Y\|^2 U (\la_0),\\
&(\la_1-\la_0) \< (\n_XJ_1)Y, U \> = - (\la_1-\la_0) \< U, (\n_XJ_1)J_1X \> = -\|X\|^2 U (\la_0) = -\|Y\|^2 U (\la_0),
\end{align*}
so \eqref{eq:nablal1} takes the form $U(\la_1+2\la_0)= 0$. Then $\la_1 + 2\la_0 = c_1$, a constant.
Similarly, $\la_2 + 2\la_0 =c_2$, $\la_3 + 2\la_0 = c_3$ for some constants $c_2, c_3$.

The manifold $M^8$ is pointwise Osserman, hence Einstein, so $\la_1 + \la_2 + \la_3 + 4\la_0 = C$,
a constant (the scalar curvature). Then $\la_0 = (c_1+c_2+c_3-C)/2$ is constant,
and the claim follows from \cite[Theorem 7.1]{GSV}.

\begin{remark} The above computation similarly works for $4$-dimensional pointwise Osserman manifolds
with $\Cliff(1)$-structure, including ``generalized complex space forms" \cite{Ols}. However, the equation
$\la_1+2\la_0=\mathrm{const}$ gives nothing new: it is just the condition of constancy
of the scalar curvature.
\end{remark}



\begin{thebibliography}{99}

\bibitem{ABS}
M.F.Atiah, R.Bott, A.Shapiro,
Clifford modules,
Topology,
\textbf{3, suppl.1} (1964), 3 -- 38.

\bibitem{Besse}
A.L.Besse,
Manifolds all of whose geodesics are closed,
Ergebnisse der Mathematik und ihrer Grenzgebiete,
\textbf{93}(1978), Springer-Verlag.

\bibitem{BG}
R.Brown, A.Gray,
Riemannian manifolds with holonomy group $Spin(9)$,
Diff. Geom. in honor of K.Yano, Kinokuniya, Tokyo (1972), 41 -- 59.

\bibitem{Chi}
Q.-S.Chi,
A curvature characterization of certain locally rank-one symmetric spaces,
J. Differ. Geom.,
\textbf{28}(1988), 187 -- 202.

\bibitem{GKV}
E.Garc{\'\i}a-R{\'\i}o, D.N. Kupeli, R.V{\'a}zguez-Lorenzo,
Osserman manifolds in Semi-Riemannian Geometry,
Lecture Notes in Mathematics,
\textbf{1777}(2002), Springer-Verlag.

\bibitem{G}
P.Gilkey,
Manifolds whose curvature operator has constant eigenvalues at the basepoint,
J. Geom. Anal.,
\textbf{4}(1994), 155 -- 158.

\bibitem{GSV}
P.Gilkey, A.Swann, L.Vanhecke,
Isoperimetric geodesic spheres and a conjecture of Osserman concerning the Jacobi
operator,
Quart. J. Math. Oxford (2),
\textbf{46}(1995), 299 -- 320.

\bibitem{H}
D.Husemoller, Fiber bundles,
(1975), Springer-Verlag.

\bibitem{LM}
H.B.Lawson, M.-L.Michelsohn, Spin geometry,
(1989), Princeton Univ. Press.

\bibitem{Lev}
J.Levine, Imbedding and immersion of real projective spaces,
Proc. Amer. Math. Soc.,
\textbf{14} (1963), 801 -- 803.

\bibitem{Nag}
M.Nagata, A remark on the unique factorization theorem,
J. Math. Soc. Japan,
\textbf{9}, (1957), 143--145.

\bibitem{Nhjm}
Y.Nikolayevsky,
Osserman manifolds and Clifford structures,
Houston J. Math.,
\textbf{29}(2003), 59--75.

\bibitem{Ndga}
Y.Nikolayevsky,
Two theorems on Osserman manifolds,
Diff. Geom. Appl.,
\textbf{18}(2003), 239 -- 253.

\bibitem{NarXiv}
Y. Nikolayevsky,
Osserman Conjecture in dimension $n \ne 8, 16$, preprint,\newline
\verb"http://au.arxiv.org/abs/math.DG/0204258"

\bibitem{Ols}
Z.Olszak,
On the existence of generalized space forms,
Israel J. Math.,
\textbf{65}(1989), 214 -- 218.

\bibitem{O}
R.Osserman,
Curvature in the eighties,
Amer. Math. Monthly,
\textbf{97}(1990), 731 -- 756.

\bibitem{Pf}
A.Pfister, Quadratic forms with applications to algebraic geometry and topology,
London Math. Soc. Lecture Notes Ser., \textbf{217}, (1995),
Cambridge Univ. Press.

\end{thebibliography}
\end{document}